\documentclass[preprint,1p]{elsarticle}
\usepackage{graphicx}
\usepackage{amssymb}
\usepackage{amsmath,amsthm}
\usepackage{amsfonts}
\usepackage{graphicx}
\usepackage{graphics}
\usepackage{color}
\usepackage{latexsym}
\usepackage{url}
\newcommand{\R}{\ensuremath{\mathbb{R}}}

\newcommand{\N}{\ensuremath{\mathbb{N}}}
\newcommand{\Z}{\ensuremath{\mathbb{Z}}}

\def\bfm#1{\protect{\makebox{\boldmath $#1$}}}
\def\x {\bfm{x}}

\def\a {\bfm{a}}

\def\u {\bfm{u}}

\def\p {\bfm{p}}

\def\bbeta {\bfm{\beta}}

\def\eep {\bfm{\varepsilon}}

\def\ep{\epsilon}

\def\wh{\widehat}

\def\sgreek#1{\protect{\makebox{\footnotesize $#1$}}}
\def\schi{\sgreek{\chi}}

\newtheorem{Theorem}{Theorem}[section]

\theoremstyle{definition}

\theoremstyle{remark}
\newtheorem{Remark}[Theorem]{Remark}
\numberwithin{equation}{section}

\begin{document}
\begin{frontmatter}
\title{A new low-cost meshfree method for two and three dimensional problems in elasticity}%
\author{Davoud Mirzaei}     
\ead{d.mirzaei@sci.ui.ac.ir}%
\address{Department of Mathematics, University of Isfahan, 81745-163 Isfahan, Iran. }


\date{\small\textsl{February 9, 2014}}
\begin{abstract}
In this paper, we continue the development of the {\em Direct Meshless Local Petrov-Galerkin (DMLPG)} method for elasto-static problems.
This method is based on the generalized moving least squares approximation. The computational efficiency is the most significant advantage of the new method in comparison with the original MLPG.
Although, the ``Petrov-Galerkin" strategy is used to build the primary local weak forms, 
the role of trial space is ignored and {\em direct approximations} for local weak forms and boundary conditions are performed to construct the final stiffness matrix. In this modification the numerical integrations are performed over polynomials instead of complicated MLS shape functions.
In this paper, DMLPG is applied for two and three dimensional problems in elasticity. Some variations of the new method are developed and their efficiencies are reported. Finally, we will conclude that DMLPG can replace the original MLPG in many situations.
\end{abstract}
\end{frontmatter}
\section{Introduction}\label{SectIntro}
The {\em Meshless Local Petrov-Galerkin (MLPG)} method has been widely employed to find the numerical solutions of
elasto-static and elasto-dynamic problems. MLPG was first introduced in \cite{atluri-zhu:1998-1}, and was first applied to elasticity in
\cite{atluri-zhu:2000-1}. Afterward, many papers were appeared for different types of mechanical problems. For example see
\cite{li-et-al:2003-1,sladek-et-al:2005-1} and the recent review paper \cite{sladek-et-al:2013-1}.
MLPG is based on {\em local weak forms} and it is known as a {\em truly} meshless method, because it uses no global background mesh
to evaluate integrals, and everything breaks down to some regular,
well-shaped and independent sub-domains. This is in contrast with methods which are based on {\em global weak forms}, such as the Element-free Galerkin (EFG) method \cite{belytschko-et-al:1994-1}, where triangulation is again required for
numerical integration. {\em But} MLPG still suffers from the cost of numerical integration. This is due to the complexity of the integrands.
In MLPG and all MLS based methods, integrations are done over complicated MLS shape functions, and this leads to high computational costs in comparison with the finite elements method (FEM), where integrands are simple and close form polynomials.
Thus, special cares should be taken in performing numerical quadratures for meshfree methods.
These challenges have been addressed in
various engineering papers \cite{beissel-belytsschko:1996-1,dolbow-belytsschko:1999-1,atluri-et-al:1999-2,carpinteri-et-al:2002-1,pecher:2006-1,mazzia-pini:2010-1} and several approaches to
implement numerical integration have been proposed in the literature. A brief review of these approaches is presented in Section 3 of \cite{babuska-et-al:2009-1}.

This is the reason why this method, and of course the other meshfree methods, have found very limited application to three-dimensional problems, which are routine applications of FEM.

A tricky modification has been applied to MLPG in \cite{mirzaei-schaback:2013-1}, in which the numerical integrations are done over low-degree polynomial basis functions rather than complicated MLS shape functions. In addition, as the shapes of the local subdomains remain unchanged, the values of integrals remain the same.
This reduces the computational costs of MLPG, significantly. In the new method, local weak forms are considered as functionals and they are {\em directly} approximated from nodal data using a generalized moving least squares (GMLS) approximation. Thus this method is called Direct MLPG (DMLPG).
Although DMLPG uses the same local forms, it is theoretically different from MLPG, because it eliminates the role of trial space. DMLPG can be considered as a generalized finite difference method (GFDM), not only in its usual strong form, but also in a weak formulation.
It is worthy to note that, by this modification we do not lose the order of convergence. This has been analytically proven in \cite{mirzaei-et-al:2012-1,mirzaei:2013-1} for different definitions of functionals, specially for the local weak forms of DMLPG.

DMLPG has been applied to the heat conduction problem in \cite{mirzaei-schaback:2013-2} and has been numerically investigated for 2D and 3D potential
problems in \cite{mazzia-et-al:2012-1}.

In this paper, the application of DMLPG is provided for elasto-static problems for the first time.
We consider both two and three dimensional problems to show the efficiency of the new method. The method can be easily extended to the other problems in elasticity.

\section{Generalized moving least squares}\label{Sect}
Generalized moving least squares (GMLS) approximation was presented in \cite{mirzaei-et-al:2012-1} in details. Here we briefly discuss this concept.
Let $\Omega$ be a bounded subset in $\R^d$, $d\in\Z_{+}$, and $X=\{x_1,x_2,\ldots,x_N\}\subset \Omega$ be a set of meshless points scattered (with certain quality) over $\Omega$. The MLS method approximates the function $u\in U$ (with certain smoothness) by its values at points $x_j$, $j=1,\ldots,N$, by
\begin{equation}\label{mls-app}
u(x)\approx \wh u(x) = \sum_{j=1}^N a_j(x) u(x_j), \quad x\in \Omega,
\end{equation}
where $a_j(x)$ are MLS shape functions obtained in such way that $\wh u(x)$ be the best approximation of $u(x)$ in polynomial subspace $\mathbb P_m(\R^d)=\mathrm{span}\{p_1,\ldots,p_Q\}$,  $Q={m+d\choose d}$, with respect to a weighted, discrete and {\em moving} $\ell^2$ norm.
The weight function governs the influence of the data points and assumed to be a function
$w:\Omega\times \Omega \to \R$ which becomes smaller
the further away its arguments are from each other. Ideally,
$w$ vanishes for arguments $x, y\in\Omega$ with $\|x - y\|_2$
greater than a certain threshold, say $\delta$. Such a behavior can be
modeled by using a translation-invariant weight function. This
means that $w$ is of the form $w(x, y) = \varphi(\|x - y\|_2/\delta)$ where
$\varphi$ is a compactly supported function supported in $[0,1]$.
If we define
\begin{equation}\label{eq-W-P}
\begin{split}
P&= P(x)=\big(p_k(x_j)\big)\in \R^{N\times Q}, \\
W&=W(x) = \mathrm{diag}\{w(x_j,x)\}\in \R^{N\times N},
\end{split}
\end{equation}
then a simple calculation
gives the \emph{shape functions}
\begin{equation}\label{mls-shapefunc}
\a(x):=[a_1(x),\ldots,a_N(x)] = \p(x) (P^TWP)^{-1}P^TW.
\end{equation}
where $\p=[p_1,\ldots,p_Q]$.
If $X_x=\{x_j: \|x-x_j\|\leqslant \delta\}$ is $\mathbb P_m(\R^d)$-unisolvent then $A(x)=P^TWP$ is positive definite \cite{wendland:2005-1} and the MLS approximation is well-defined at sample point $x$.
Of course if $\|x-x_j\|\geqslant \delta$ then $a_j(x)=0$. Thus, in programming we can only form $P$ and $W$ for
active points $X_x$ instead of $X$.
Derivatives of $u$ are usually approximated by derivatives of $\wh u$,
\begin{equation}\label{mls-standard-deriv}
D^\alpha u(x)\approx D^\alpha \wh u(x) = \sum_{j=1}^N D^\alpha a_j(x) u(x_j), \quad x\in \Omega, \quad \alpha=(\alpha_1,\ldots,\alpha_d)\in\N_0^d.
\end{equation}
These derivatives are sometimes called {\em standard} or {\em full} derivatives. Details are in \cite{lancaster-salkauskas:1981-1,belytschko-et-al:1996-1,mirzaei:2015-1} and any other text containing the application of MLS approximation.

The GMLS approximation can be introduced as below. Suppose that $\lambda$ is a linear functional from the dual space $U^*$.
The problem is the recovery of $\lambda(u)$ from nodal values $u(x_1),\ldots,u(x_N)$.
The functional $\lambda$ can, for
instance, describe point evaluations of $u$, its derivatives up to order $m$, and the weak formulations which involve $u$ or a
derivative against some test function. The approximation $\wh \lambda(u)$ of
$\lambda (u)$ should be a linear function of the data $u(x_j)$, i.e., it should have the form
\begin{equation}\label{gmls-app}
\lambda(u)\approx \wh \lambda (u)=\sum_{j=1}^Na_j(\lambda)u(x_j),
\end{equation}
where $a_j(\lambda)$ are shape functions associated to the functional $\lambda$. If $\lambda$ is chosen to be the {\em point evaluation} functional
$\delta_x$, where $\delta_x(u):=u(x)$, then the classical MLS approximation \eqref{mls-app} is obtained.
If we assume $\lambda$ is finally evaluated at sample point $x$, then the same weight function $w(x,y)$ as in the classical MLS can be used which is independent of the choice of $\lambda$. Using this assumption,
analogous to \eqref{mls-shapefunc}, \cite{mirzaei-et-al:2012-1} proves,
\begin{equation}\label{gmls-shapefunc}
\a(\lambda):=[a_1(\lambda),\ldots,a_N(\lambda)] = \lambda(\p) (P^TWP)^{-1}P^TW,
\end{equation}
where $\lambda(\p)=[\lambda(p_1), \ldots, \lambda(p_Q)]$.
In fact, we have a {\em direct} approximation for $\lambda(u)$ from nodal values $u(x_1),\ldots,u(x_N)$, without any detour via classical MLS shape functions.
One can see, $\lambda$ acts only on polynomial basis functions. This is the central idea in this GMLS approximation which finally speeds up our numerical algorithms. If $\lambda$ contains derivatives of $u$, \eqref{gmls-shapefunc} shows that derivatives of weight functions are not required. This paves the way to generalize the forthcoming schemes for discontinuous problems.

In particular, if $\lambda(u) = D^\alpha(u)$ then derivatives of $u$ are recovered. They are different from the standard derivatives \eqref{mls-standard-deriv}, and in meshless literature they are called {\em diffuse} or {\em uncertain} derivatives. But \cite{mirzaei-et-al:2012-1} and \cite{mirzaei:2013-1} prove the optimal rate of convergence for them toward the exact derivatives, and thus there is nothing diffuse or uncertain about them. As suggested in
\cite{mirzaei-et-al:2012-1}, they can be called {\em GMLS derivative} approximations.

In the next sections, we deliberately choose $\lambda$ in such way that MLPG methods speed up, significantly.

The GMLS approximation of this section is different from one presented in \cite{atluri-et-al:1999-1}. In that paper a Hermite-type MLS approximation has been used to solve the forth order problems of thin beams.
Here we approximate the general functional $\lambda(u)$ from values
$u(x_1),\ldots, u(x_N)$, where information of $D^\alpha u$ is not required.
In more general situation, the GMLS approximation of \cite{atluri-et-al:1999-1} can be written as
\begin{equation*}
u(x)\approx \wh u(x)=\sum_{k=1}^K\sum_{j=1}^N a_{k,j}(x) \mu_{k,j}(u),
\end{equation*}
where $\mu_{k,j}$ are linear functionals from $U^*$ and should be chosen properly to ensure the solvability of the problem.

In a more and more general situation, both these generalizations can be used simultaneously
\begin{equation*}
\lambda (u)\approx \wh \lambda (u)=\sum_{k=1}^K\sum_{j=1}^N a_{k,j}(\lambda) \mu_{k,j}(u).
\end{equation*}
So far, there is no rigorous error analysis for such generalized approximation, even when $\lambda$ and $\mu_{k,j}$ are some special functionals. Throughout, we leave the above recent formulations and focus on GMLS approximation \eqref{gmls-app} together with \eqref{gmls-shapefunc}.

\section{Local weak forms of the elasticity problem}\label{Sect}
Let $\Omega\subset\R^d$ (usually $d=2,3$) be a bounded domain with boundary $\Gamma$. From here on, integers $i$ and $j$ are assumed to vary from $1$ to $d$. Consider the following $d$-dimensional elasto-static
problem
\begin{equation}\label{elasto-static-eq}
\sigma_{ij,j}+b_i = 0,\quad \mbox{in }\Omega \,
\end{equation}
where $\sigma_{ij}$ is the stress tensor, which corresponds to the displacement field $u_i$, and
$b_i$ is the body force.
The corresponding boundary conditions are given by
\begin{align}
u_i &= \overline u_i,\quad \mbox{on }\Gamma_u,\\
t_i &= \sigma_{ij}n_j = \overline t_i,\quad \mbox{on }\Gamma_t,
\end{align}
where $\overline u_i$ and $\overline t_i$ are the prescribed displacement and traction on
the boundaries $\Gamma_u$ and $\Gamma_t$, respectively. $n$ is the unit outward
normal to the boundary $\Gamma$.

Many numerical methods such as FEM, FVM, BEM, EFG, etc. are based on a global weak form of
\eqref{elasto-static-eq} over entire $\Omega$, which can be derived using the integration by parts. However, the MLPG method starts
from weak forms over sub-domains $\Omega_k$
inside the global domain $\Omega$. Sub-domains usually cover the entire domain $\Omega$ and they have simple geometries in order to do the numerical integrations as easily as possible.

Let
$X=\{x_1,x_2,\ldots,x_N\}\subset \Omega$ be a set of scattered meshless points, where some points are located on the boundary $\Gamma$ to enforce the boundary conditions.
In this work, spherical (circular in 2D) subdomains $\Omega_k=B(x_k,r_k)\cap \Omega$ with radius $r_k$ centered at $x_k$, and cubical (rectangular in 2D) subdomains
$\Omega_k=C(x_k,s_k)\cap\Omega$ with side-length $s_k$ centered at $x_k$ are employed.
Of course, for boundary points, $\partial \Omega_k$ intersects with the global boundary $\Gamma$.
A local weak form of the equilibrium
equation over $\Omega_k$ is written as
\begin{equation}\label{lwf1}
\int_{\Omega_k}(\sigma_{ij,j}+b_i)v_i \, d\Omega = 0,
\end{equation}
where $v_i$ are appropriate test functions. We do not introduce Lagrange multiplier or penalty parameter in the weak form, because
in our numerical method the essential boundary conditions are imposed in a suitable collocation form. Thus we assume $x_k$ is located either inside $\Omega$ or on $\Gamma_t$ where the tractions are prescribed. Using $\sigma_{ij,j}v_i=(\sigma_{ij}v_i)_{,j}-\sigma_{ij}v_{i,j}$ and
the Divergence Theorem, from \eqref{lwf1} we have
\begin{equation}\label{lwf2}
\int_{\partial\Omega_k}\sigma_{ij}n_jv_i\, d\Gamma -\int_{\Omega_k}\sigma_{ij}v_{i,j} \, d\Omega =
\int_{\Omega_k}b_iv_{i} \, d\Omega,
\end{equation}
where $n$ is the outward unit normal to the boundary $\partial\Omega_k$. Imposing the natural boundary conditions $\sigma_{ij}n_j=\overline t_i$
on $\partial\Omega_k\cap \Gamma_t$, we have
\begin{equation}\label{lwf3}
\int_{\partial\Omega_k\setminus \Gamma_t}\sigma_{ij}n_jv_i\, d\Gamma -\int_{\Omega_k}\sigma_{ij}v_{i,j} \, d\Omega =
 \int_{\Omega_k}b_iv_{i} \, d\Omega -\int_{\partial\Omega_k\cap \Gamma_t}\overline t_i v_i\, d\Gamma.
\end{equation}
In Petrov-Galerkin methods, the
trial functions and the test functions come from different spaces. Thus there will be many choices for test functions $v_i$, and this leads to
a list of MLPG methods labeled from 1 to 6. But this may cause some difficulties in mathematical analysis.
Up until here, the new procedure is identical to the classical MLPG method.
In the next section we pave the way of going from MLPG to DMLPG using the concept of GMLS approximation.
\section{DMLPG formulation}\label{Sect}
Although, DMLPG uses the same local weak forms obtained from a Petrov-Galerkin formulation, it is mathematically different
from MLPG because direct approximations for local weak forms are provided to rule out the action of trial space.

Using the same labels as in MLPG, here we discuss DMLPG1 and 5 and leave the others for a new research.
Note that
there are some difficulties to develop DMLPG3 and 6 because they are based on a Galerkin formulation \cite{mirzaei-schaback:2013-1,mirzaei-schaback:2013-2}.

We use the same scheme to impose the essential boundary conditions in all types of DMLPG. The MLS collocation method is applied at
points located on $\Gamma_u$,
\begin{equation}\label{essBCimpose}
\sum_{\ell=1}^N a_\ell(x_k)u_i(x_\ell) = \overline u_i(x_k), \quad x_k\in \Gamma_u.
\end{equation}
In fact, the functional $\lambda$ in GMLS is taken to be $\delta_{x_k}$, the point evaluation functionals at $x_k$.
In the following subsections, we consider the local weak forms around the points located either inside $\Omega$ or over Neumann parts of the boundary $\Gamma$.
\subsection{DMLPG1}
Let $\u=[u_1,\ldots,u_d]^T$. If test functions $v_i$ are chosen such that they all vanish over $\partial\Omega_k\setminus\Gamma_t$, then
the first integral in \eqref{lwf3} vanishes and if we define
\begin{equation}\label{lamdak1}
\begin{split}
\lambda_{k}^{(i)}(\u) :=& -\int_{\Omega_k}\sigma_{ij}v_{i,j} \, d\Omega , \\
\beta_{k}^{(i)} :=&\int_{\Omega_k}b_iv_{i}\, d\Omega -\int_{\partial\Omega_k\cap \Gamma_t}\overline t_i v_i\, d\Gamma,
\end{split}
\quad\quad  x_k\in \mathrm{int}(\Omega)\cup \Gamma_t,
\end{equation}
then\eqref{lwf2} becomes
$$
\lambda_{k}^{(i)}(\u) = \beta_{k}^{(i)},\quad x_k\in \mathrm{int}(\Omega)\cup \Gamma_t.
$$
Now, the GMLS can be applied to approximate the above functionals.
To simplify the notation, let
$$
\bbeta_k =\begin{bmatrix}\beta_k^{(1)}\\ \vdots\\ \beta_k^{(d)}\end{bmatrix},
\quad
\u= \begin{bmatrix}u_1\\ \vdots\\ u_d
\end{bmatrix},\quad A_{k\ell}= \begin{bmatrix} a_{k\ell}^{(11)}& \cdots &a_{k\ell}^{(1d)}\\ \vdots & \ddots &\vdots \\ a_{k\ell}^{(d1)} &\cdots &a_{k\ell}^{(dd)}\end{bmatrix},
$$
where $A=(A_{k\ell})$ is introduced as a block matrix for reserving the acts of GMLS functions. Blocks of $A$ are not diagonal, because $\lambda_{k}^{(i)}(\u)$ depends not only on $u_i$ (for a specified $i$) but also on all $u_i$ for $i=1,\ldots,d$.
The GMLS approximation can be used to write
\begin{equation}\label{gmlsOffunctionals-block}
\lambda_{k}(\u)\approx \wh{\lambda}_{k}(\u)=\sum_{\ell=1}^N A_{k\ell}\u(x_\ell).
\end{equation}
According to \eqref{gmls-shapefunc}, if $A_{k,:}$ represents the $k$-th block row of $A$, then
\begin{equation}
A_{k,:} = \lambda_{k}(\p) \Phi \in \R^{d\times dN},
\end{equation}
where $\Phi\in \R^{dQ\times dN}$ is a block matrix obtained from $\phi := (P^TWP)^{-1}WP^T \in \R^{Q\times N}$ by
$$
\Phi_{ij} = \begin{bmatrix} \phi_{ij}& &0\\ &\ddots & \\ 0 & &\phi_{ij}\end{bmatrix}\in \R^{d\times d}.
$$
Matrices $P$ and $W$ are defined in \eqref{eq-W-P}, and $\p$ is defined by
$$
\p = \begin{bmatrix}  p_1(x)&p_2(x)&\cdots &p_Q(x)\\
\vdots & \vdots &&\vdots \\p_1(x)&p_2(x)&\cdots &p_Q(x)\\
\end{bmatrix}\in\R^{d\times Q}.
$$
Thus we have
\begin{equation}\label{lambda_approx1}
\lambda_{k}(\p)
= -\Big[\underbrace{\int_{\Omega_k}\eep_v D P_1(x) d\Omega}_{\in\R^{d\times d}} ,\underbrace{ \int_{\Omega_k}\eep_v D P_2(x) d\Omega}_{\in\R^{d\times d}} ,\ldots, \underbrace{\int_{\Omega_k}\eep_v D P_Q(x) d\Omega}_{\in\R^{d\times d}}
\Big]\in \R^{d\times dQ},
\end{equation}
where for a two dimensional problem ($d=2$) of
isotropic material, the stress-strain matrix $D$ is defined by
$$
D=\frac{\overline E}{1-\overline \nu^2}\begin{bmatrix} 1& \overline \nu &0\\ \overline \nu &1 &0 \\0& 0&(1-\overline \nu)/2 \end{bmatrix},
$$
where
$$
\overline E=\begin{cases}E & \mbox{for plane stress }\\\frac{E}{1-\nu^2} &\mbox{for plane strain }  \end{cases}\quad
\overline \nu=\begin{cases}\nu & \mbox{for plane stress }\\\frac{\nu}{1-\nu} &\mbox{for plane strain }  \end{cases},
$$
in which $E$ and $\nu$ are Youngs modulus and Poissons ratio, respectively. The strain matrix for test functions $v_i$ is
$$
\eep_{v} = \begin{bmatrix} v_{1,1}& 0 & v_{1,2}\\ 0 & v_{2,2}  & v_{2,1} \end{bmatrix},
$$
and
$$
P_n(x) = \begin{bmatrix}p_{n,1}(x)&0\\0& p_{n,2}(x)\\p_{n,2}(x)&{p_{n,1}}(x)\end{bmatrix},\quad n=1,2,\ldots,Q.
$$
For the elasticity problem of isotropic material in 3D (i.e. $d=3$), we have $D=\begin{bmatrix} D_1&0\\0&D_2 \end{bmatrix}\in \R^{6\times6}$ where
$$
D_1=\frac{E}{(1-2\nu)(1+\nu)}\begin{bmatrix} 1-\nu&\nu & \nu \\ \nu &1-\nu &\nu\\
\nu&\nu&1-\nu \end{bmatrix}, \quad D_2 = \frac{E}{2(1+\nu)}\begin{bmatrix} 1&0&0\\0&1&0\\ 0&0&1\end{bmatrix}.
$$
In addition, the strain matrix of test functions is
$$
\eep_{v} = \begin{bmatrix} v_{1,1}& 0 & 0 &0 &v_{1,3}& v_{1,2}\\ 0 & v_{2,2} &0 &v_{2,3}&0 & v_{2,1}
\\ 0&0& v_{3,3}& v_{3,2}& v_{3,1}&0 \end{bmatrix},
$$
and finally
$$
P_n(x) = \begin{bmatrix}p_{n,1}(x)&0 & 0 \\0& p_{n,2}(x)& 0 \\0 & 0 & p_{n,3}(x) \\
0 & p_{n,3}(x)& p_{n,2}(x) \\ p_{n,3}(x)& 0& p_{n,1}(x) \\  p_{n,2}(x) & p_{n,1}(x) &0\end{bmatrix},\quad n=1,2,\ldots,Q.
$$
For simplicity we choose $v_1=\cdots=v_d=:v$ in the following numerical algorithms.
To set up the final linear system, we first assume
$$\u = [u_1(x_1), \ldots,u_d(x_1),u_1(x_2),\ldots,u_d(x_2),\ldots,u_1(x_N),\ldots,u_d(x_N)]^T\in \R^{dN\times 1}.$$
Without loss of generality, let the first $N_b$ meshless points are located on $\Gamma_u$. The boundary matrix $B\in\R^{dN_b\times dN}$ corresponding to the essential boundary conditions is a block matrix in which
$$
B_{k\ell}=\begin{bmatrix}
a_\ell(x_k) & & 0\\ &\ddots & \\ 0 & & a_\ell(x_k)
\end{bmatrix}_{d\times d},
$$
where $a_\ell(x_k)$ are the values of GMLS shape functions defined in \eqref{essBCimpose}. Finally, if we set
$$K=\begin{bmatrix}B\\ A\end{bmatrix}_{dN\times dN}, \quad  R = \begin{bmatrix}\overline \u(x_1) & \cdots & \overline \u(x_{Nb})& \bbeta_{N_b+1} &\cdots &\bbeta_{N} \end{bmatrix}^T_{dN\times 1},$$
then we have the final system of linear equations
\begin{equation}\label{final-linsys}
K\u = R.
\end{equation}
Sometimes, in a boundary point $x_k$, tractions $t_i$, $i\in\{i_1,\ldots,i_s\}\subset \{1,2,\ldots,d\}$, and displacements $u_i$, $i\in \{1,2,\ldots,d\}\setminus \{i_1,\ldots,i_s\}$, are prescribed.
In this case, since
the essential boundary conditions are applied using the collocation method, in the $k$-th block row of $A$, rows $i_1,\ldots, i_s$ should be replaced
by corresponding MLS shape function vectors, say $\a^{(i_m)}_{k}$, $1\leqslant m\leqslant s$, of size $dN$. These vectors are introduced as follows: first we define $\a^{(i_m)}_{k}$ as zero $dN$-vectors. Then vector components $a_1(x_k),a_2(x_k),\ldots,a_N(x_k)$ of MLS shape function are substituted in to the component indices
$i_m, i_m+d, \ldots,i_m+(N-1)d$ of $\a^{(i_m)}_{k}$. Of course the corresponding right-hand sides should form by known boundary values $\overline u_{i_m}$ instead of $\beta_{k}^{(i_m)}$.
\begin{Remark}\label{rem_1}
In DMLPG process, integrations are only appeared in \eqref{lambda_approx1}, where they are done over polynomials rather than MLS shape functions.
This is the main idea behind the DMLPG approach. In fact, DMLPG shifts the numerical
integration into the MLS itself, rather than into an outside loop over calls to MLS
routines. Moreover, if the shifted polynomial basis functions are used and if the same weight function $v$ is employed for all local sub-domains 
then $\lambda_k(\p)=\lambda_j(\p)$ provided that $\Omega_k = \Omega_j$. For example, for all interior test points only one integral should be computed if all interior local sub-domains have the same shape. 
Therefore DMLPG is extremely faster than the original MLPG.
\end{Remark}
Moreover, in some situations, we can get the exact numerical integrations with a few number of Gaussian points. For example, if cubical subdomains with polynomial test function $v$  are used in DMLPG1, the integrands are $d$-variate polynomials
of degree $(m-1)\times(n-1)$, where $n$ is the degree of the polynomial test function.
Thus a $\left\lceil\frac{(m-1)(n-1)+1}{2}\right\rceil$--point Gauss quadrature in each axis is enough for doing the exact numerical integration.
As a polynomial test function on the square or cube for DMLPG1 with $n=2$, we can use
\begin{equation}\label{cube-testfunc}
\displaystyle
v=v(x;x_k)=
\begin{cases}
\prod_{i=1}^d\left(1-\frac{4}{s_k^2}
(\schi_i-\schi_{ki})^2\right) ,& x\in C(x_k,s_k), \\
0,&  \mbox{otherwise}
\end{cases}
\end{equation}
where $x=(\schi_1, \ldots ,\schi_d)$ and
$x_k=(\schi_{k1},\ldots,\schi_{kd})$.
Note that, we should be careful for points located on the curved parts of the boundary.

\subsection{DMLPG5}
If $v=v_i\equiv1$ are chosen over $\Omega_k$, then the second integral in
\eqref{lwf3} vanishes, and by defining
\begin{equation}\label{lamdak2}
\lambda_{k}^{(i)}(\u) := \int_{\partial\Omega_k\setminus\Gamma_t}\sigma_{ij}n_j\, d\Gamma,\quad
\beta_{k}^{(i)} :=\int_{\Omega_k}b_i d\Omega-\int_{\partial\Omega_k\cap \Gamma_t}\overline t_i \, d\Gamma,  \quad x_k\in \mathrm{int}(\Omega)\cup \Gamma_t,
\end{equation}
we have
$$
\lambda_{k}^{(i)}(\u) = \beta_{k}^{(i)}.
$$
As before, we apply the GMLS to find direct approximations for functionals $\lambda_{k}^{(i)}$. Equations are the same as those where obtained for DMLPG1, except \eqref{lambda_approx1} which should be replaced by
\begin{equation}\label{lambda_approx5}
\lambda_{k}(\p)
= \Big[\int_{\partial \Omega_k\setminus \Gamma_t} \mathcal N D P_1(x) d\Gamma , \int_{\partial \Omega_k\setminus \Gamma_t}\mathcal N D P_2(x) d\Gamma ,\ldots, \int_{\partial \Omega_k\setminus \Gamma_t} \mathcal N D P_Q(x) d\Gamma
\Big]\in \R^{d\times dQ},
\end{equation}
where $ \mathcal  N$ is reserved for matrix of components of normal vector, which is  defined for the two dimensional problem by
$$
\mathcal N = \begin{bmatrix} n_{1}& 0 & n_{2}\\ 0 & n_{2}  & n_{1} \end{bmatrix},
$$
and for the three dimensional problem by
$$
\mathcal N = \begin{bmatrix} n_{1}& 0 & 0 &0 &n_{3}& n_{2}\\ 0 & n_{2} &0 &n_{3}&0 & n_{1}
\\ 0&0& n_{3}& n_{2}& n_{1}&0 \end{bmatrix}.
$$
We note that, DMLPG5 has the features mentioned in Remark \ref{rem_1} for DMLPG1. In addition,
one can see the integrals in \eqref{lambda_approx5} are all boundary integrals. Thus DMLPG5 is slightly faster.
Again if cubes are used as subdomains, a $\left\lceil\frac{m}{2}\right\rceil$--point
Gauss quadrature in each axis gives the exact solution for local boundary
integrals.

In the following section, some numerical experiments in two and three dimensional elasticity are presented to show the efficiencies of the new methods.

\section{Numerical results}
The following compactly supported Gaussian weight function
is used
$$
\displaystyle
w(x,y)=\varphi(r)=\frac{\exp(-(\epsilon r)^2)-\exp(-\epsilon^2)}{1-\exp(-\epsilon^2)}, \quad 0\leqslant r=\frac{\|x-y\|_2}{\delta}\leqslant 1,
$$
where the shape parameter $\epsilon$ is taken to be $4$ in this paper. Here $\delta=\delta(x)$ is the radius of circular (in 2D) or spherical (in 3D) support of weight function $w$ at point $x$ in question. $\delta$ should be large enough to ensure the regularity of the moment matrix $P^TWP$ in MLS/GMLS approximation. Thus $\delta$ is proportional to $h$ (mesh-size) and $m$, say $\delta=c mh$. If we have a varying-density data point,
the support size $\delta$ can vary from point to point in $\Omega$.
The polynomial degree $m=2$ and both spherical and cubical subdomains are used. For spheres, the above Gaussian weight function
with $\delta$ being replaced by the
radius $r_k$ of the local domain $\Omega_k$, is used as a test function, while for cubes, \eqref{cube-testfunc} is applied.

Displacement and strain energy relative errors will be presented in the following numerical examples. They are defined as
$$
r_{u}=\frac{\|\u^{\mathrm{exact}}-\u^{\mathrm{numerical}} \|}{\|\u^{\mathrm{exact}}\|},\quad
r_{\ep}=\frac{\|\eep^{\mathrm{exact}}-\eep^{\mathrm{numerical}} \|}{\|\eep^{\mathrm{exact}}\|},
$$
where $\|\cdot\|$ is a discrete 2-norm on a very fine mesh point in the domain $\Omega$.

All routines are written using \textsc{Matlab}$^\copyright$
 and run on a Pentium 4 PC with 8.00 GB of Memory and a 7--core 2.4 GHz CPU.

Here we should note that, the following examples may be handled by the classical techniques such as FEM and BEM with available subroutines.
However the aim of this paper is to introduce the DMLPG for elasticity problems, whereas considering the abilities of the method for
more complicated problems, such as those with discontinuity and cracks, etc., remains for new researches.

\subsection{Cantilever beam}
\begin{figure}[hbt]
\begin{center}
\includegraphics[width=10cm]{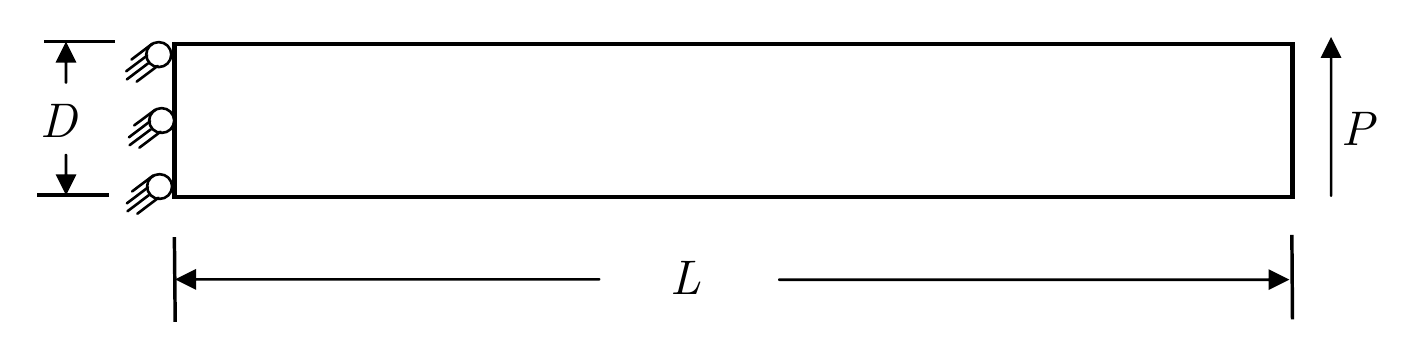}
\caption{\small{A cantilever beam}}\label{figure1}
\end{center}
\end{figure}

\begin{figure}[hbt]
\begin{center}
\includegraphics[width=8cm]{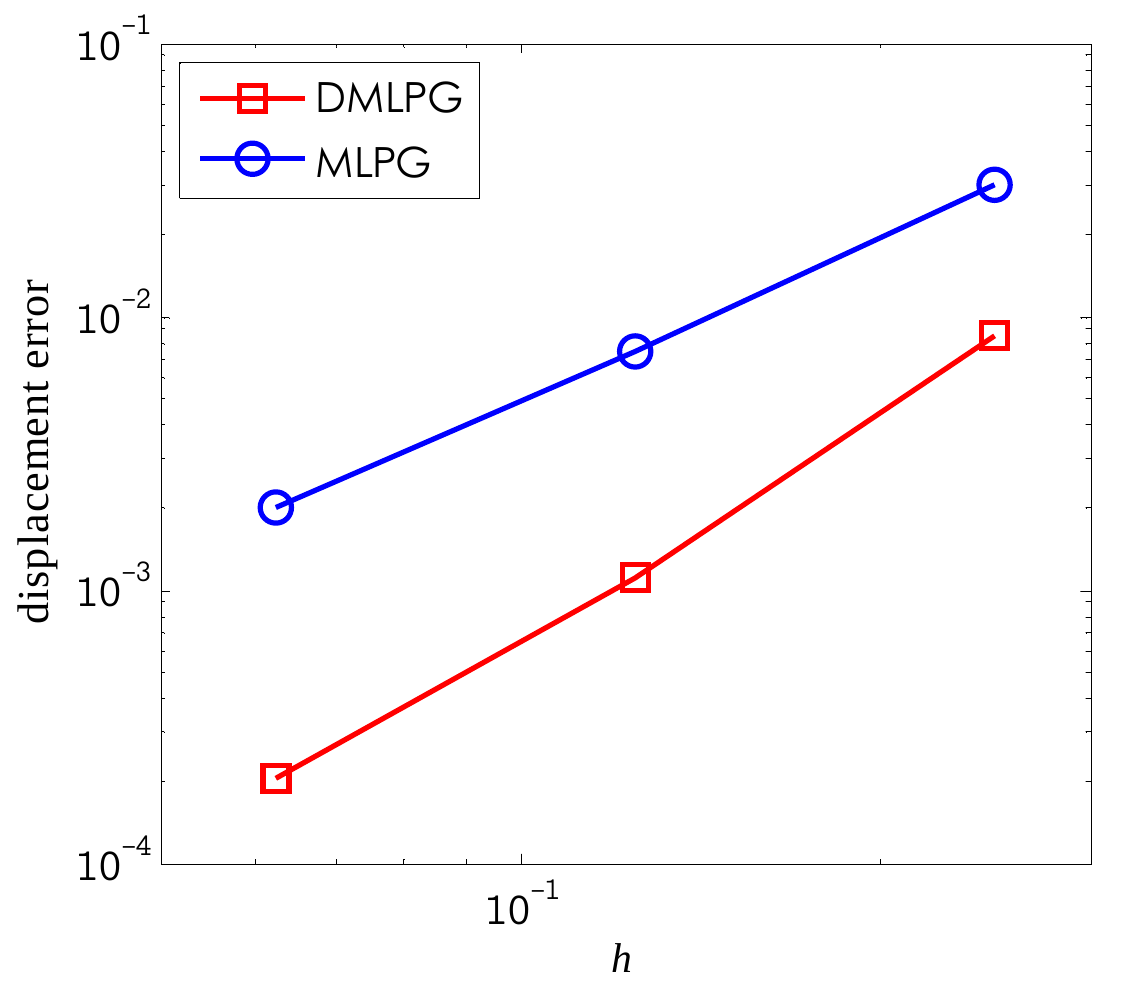}
\includegraphics[width=8cm]{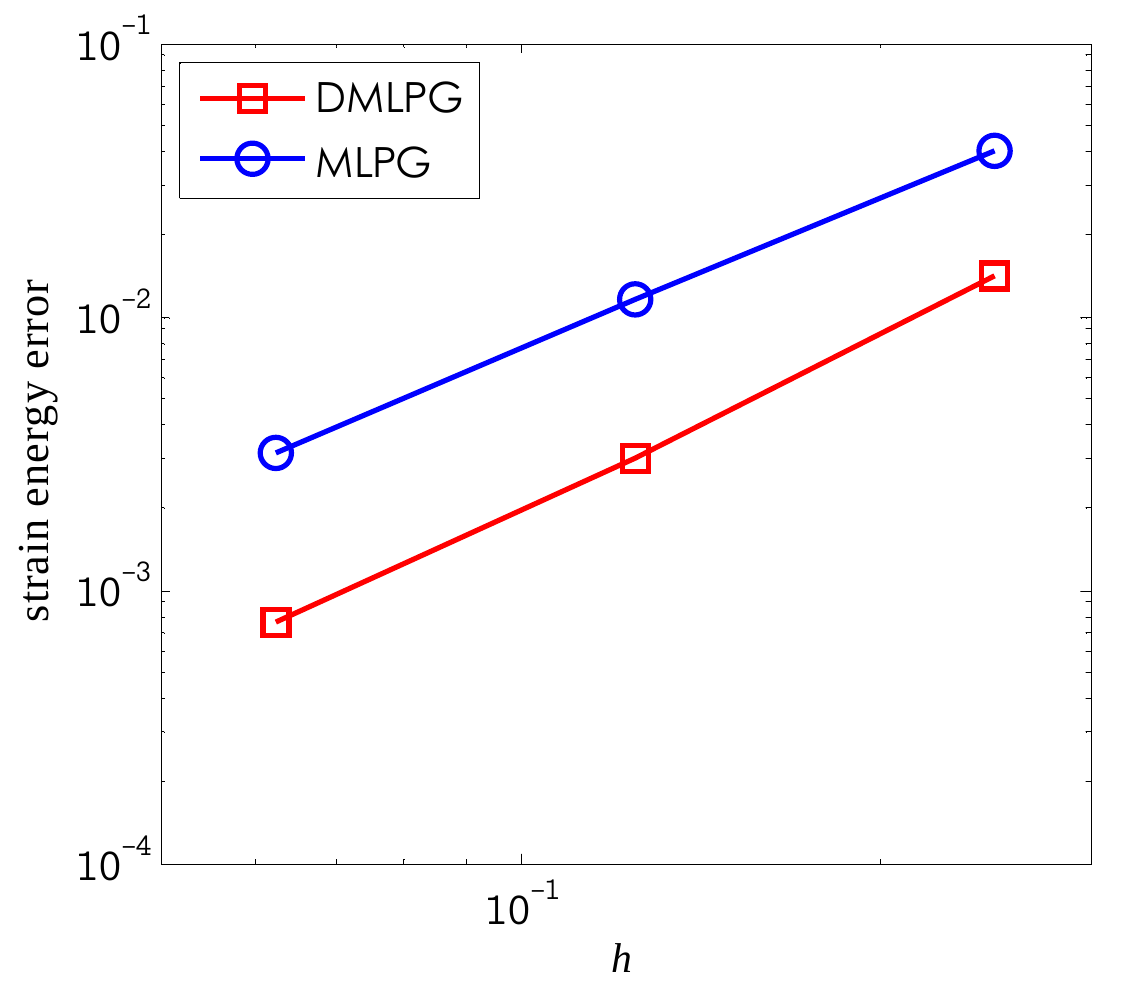}
\caption{\small{Relative displacement and strain errors for beam, rectangular subdomains }}\label{figure2}
\end{center}
\end{figure}

\begin{figure}[hbt]
\begin{center}
\includegraphics[width=8cm]{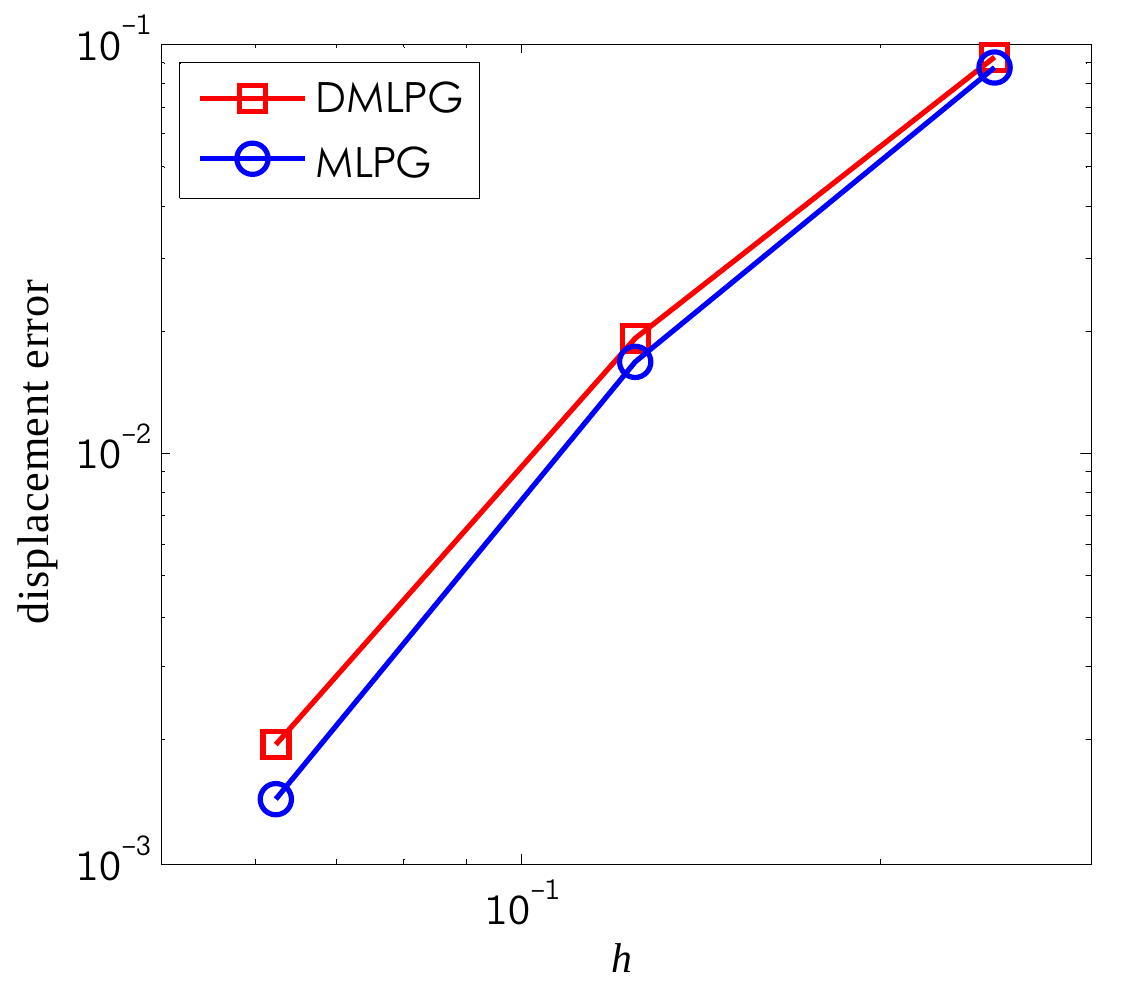}
\includegraphics[width=8cm]{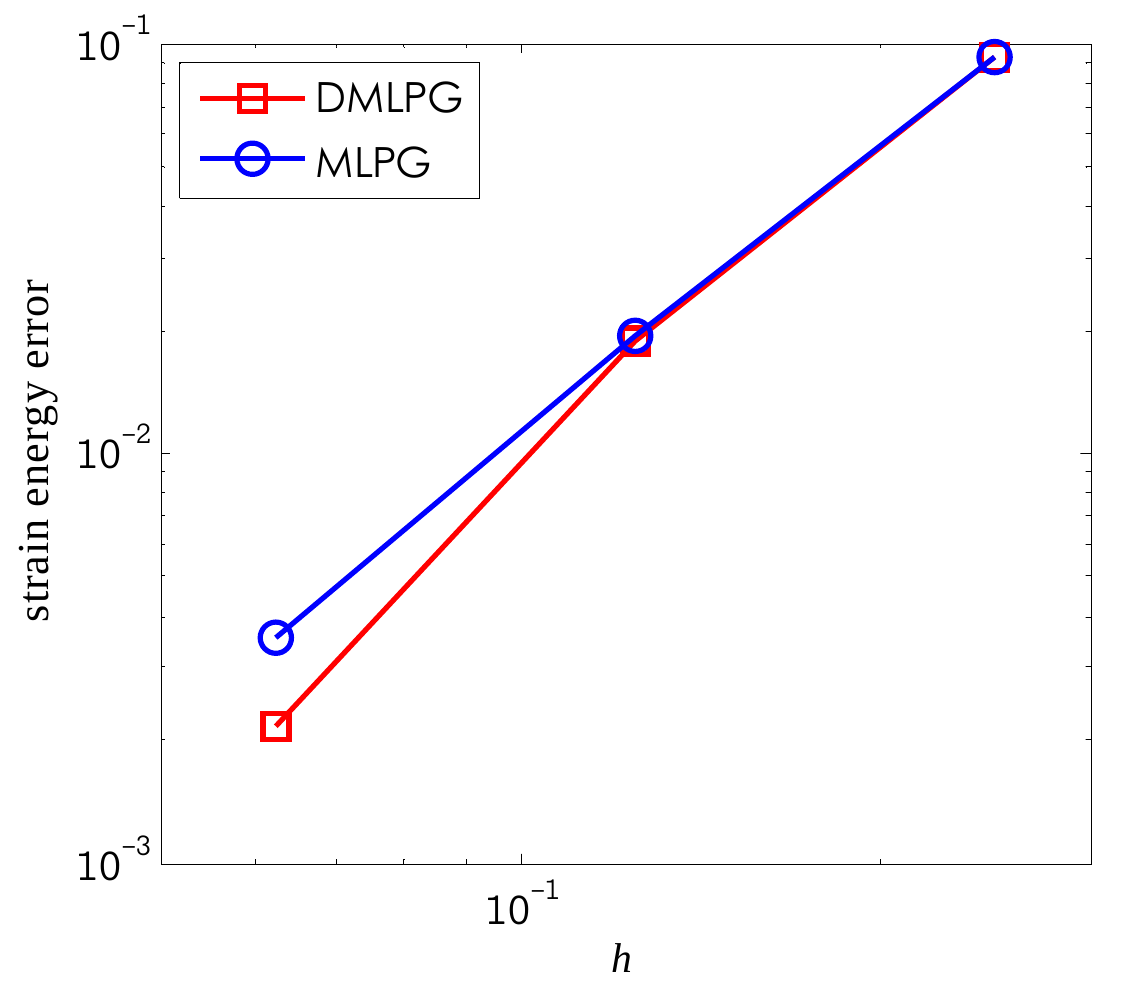}
\caption{\small{Relative displacement and strain errors for beam, circular subdomains }}\label{figure3}
\end{center}
\end{figure}

As a benchmark problem in 2D elasticity, a cantilever beam loaded by a tangential traction on the free
end, as shown in Fig. \ref{figure1}, is now considered. The exact solution of this problem is given in Timoshenko and Goodier \cite{timoshenko-goodier:1970-1} as follows:
\begin{align*}
u_1&=-\frac{P}{6\overline E I}\left(\schi_2-\frac{D}{2}\right)\big(3\schi_1(2L-\schi_1)+(2+\overline \nu)\schi_2(\schi_2-D)\big),\\
u_2&=\frac{P}{6\overline E I}\left[\schi_1^2(3L-\schi_1)+3\overline \nu(L-\schi_1)\left(\schi_2-\frac{D}{2}\right)^2+\frac{4+5\overline \nu }{4}D^2\schi_1\right],
\end{align*}
where $I=D^3/12$ and $x=(\schi_1,\schi_2)\in\R^2$. The corresponding exact stresses are
\begin{align*}
\sigma_{11}&=-\frac{P}{I}(L-\schi_1)\left(\schi_2-\frac{D}{2}\right),\\
\sigma_{22}&=0,\\
\sigma_{12}&=-\frac{P\schi_2}{2I}\left(\schi_2-D\right).
\end{align*}
Both MLPG1 and DMLPG1 are applied with $L=8$, $D=1$, $P=1$, $E=1$, $\nu=0.25$ for the plane stress case. The uniform mesh sizes
$(33\times 5)$, $(65\times 9)$ and $(129\times 17)$ are used to detect the rates of convergence and computational costs of both techniques.
Circular domains with radius $r_k=0.7h$, and rectangular domains with height-length $h\times h$ are employed as sub-domains $\Omega_k$ for all $k$. As pointed before, for $m=2$ a 2-point Gaussian quadrature in each axis is enough to get the exact numerical integrations over squares in DMLPG. But
10-point quadrature in each axis is used for circles ($r$ and $\theta$ directions) in both methods and for squares in MLPG. The sufficiently large number of Gaussian points
should be used to get the high accuracy for integration against MLS shape functions in MLPG. However, DMLPG works properly with fewer integration points, because there is no shape function incorporated in integrands. Here, to make the comparisons regarding the computational costs, we use the same number of Gaussian points for both methods in circular subdomains. Results are presented in Figs. \ref{figure2} and \ref{figure3} to compare the accuracy of numerical displacements, numerical strains
in MLPG1 and DMLPG1 for square and circle sub-domains. The rates seem to be the same, although, the results of DMLPG with
squares are more accurate. This is expectable, because in this case the integrals are computed exactly.

\begin{figure}[hbt]
\begin{center}
\includegraphics[width=8cm]{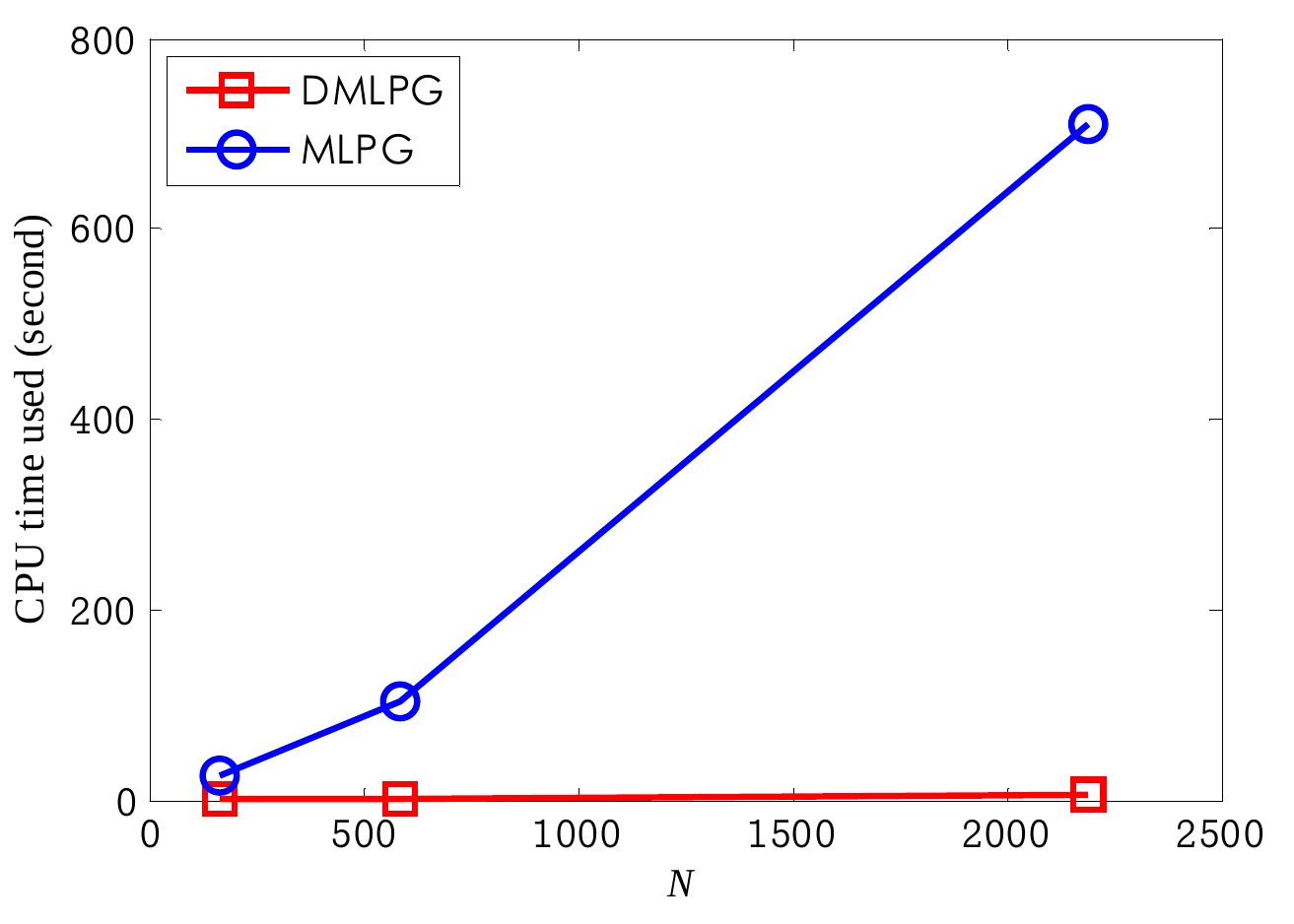}
\includegraphics[width=8cm]{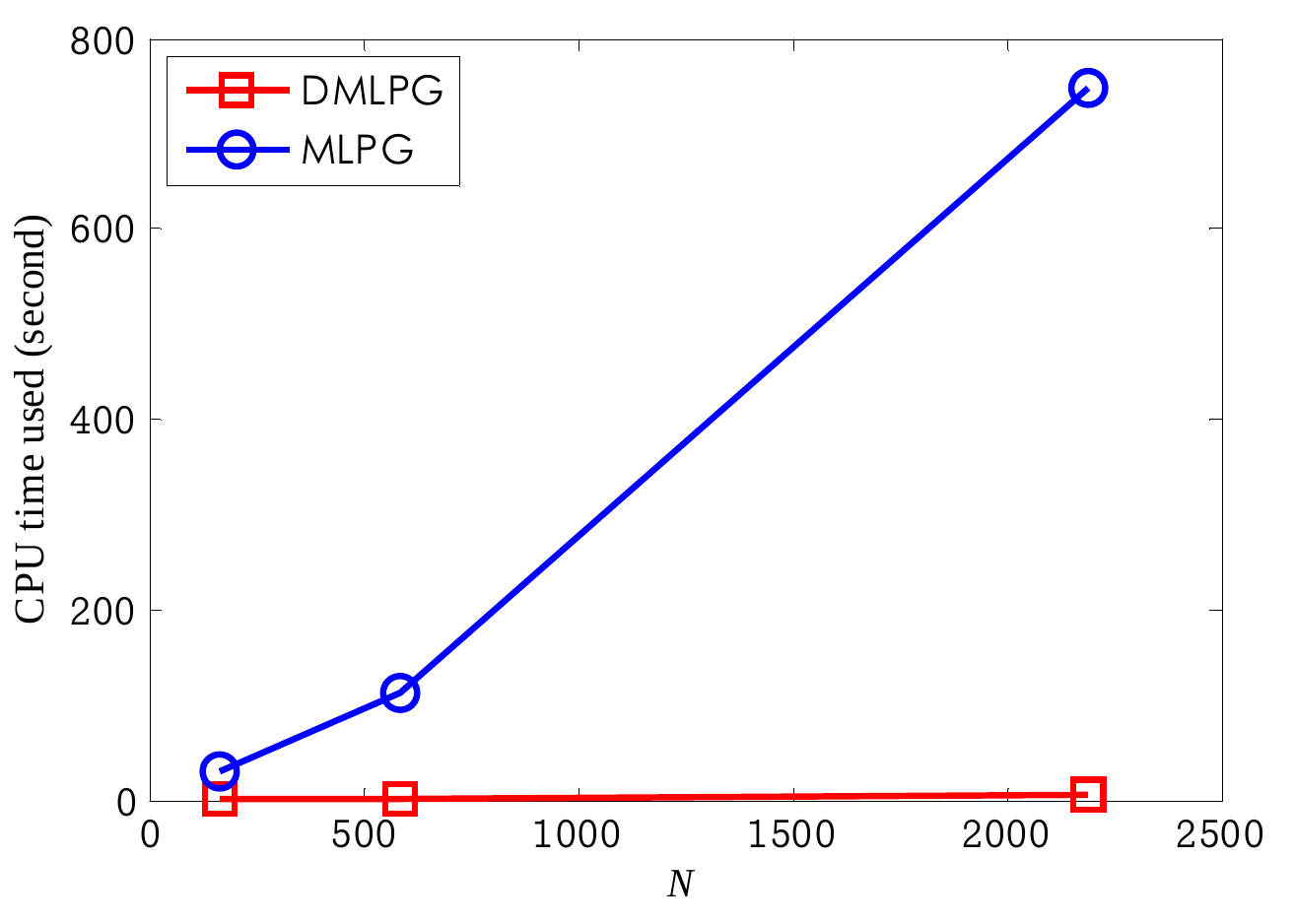}
\caption{\small{Computational costs for beam, rectangular (up) and circular (down) subdomains }}\label{figure4}
\end{center}
\end{figure}

As discussed before, DMLPG is superior to MLPG in computational efficiency. To confirm this numerically, the CPU times used are
compared in Fig.
\ref{figure4} for square and circle subdomains.

\begin{figure}[hbt]
\begin{center}
\includegraphics[width=8cm]{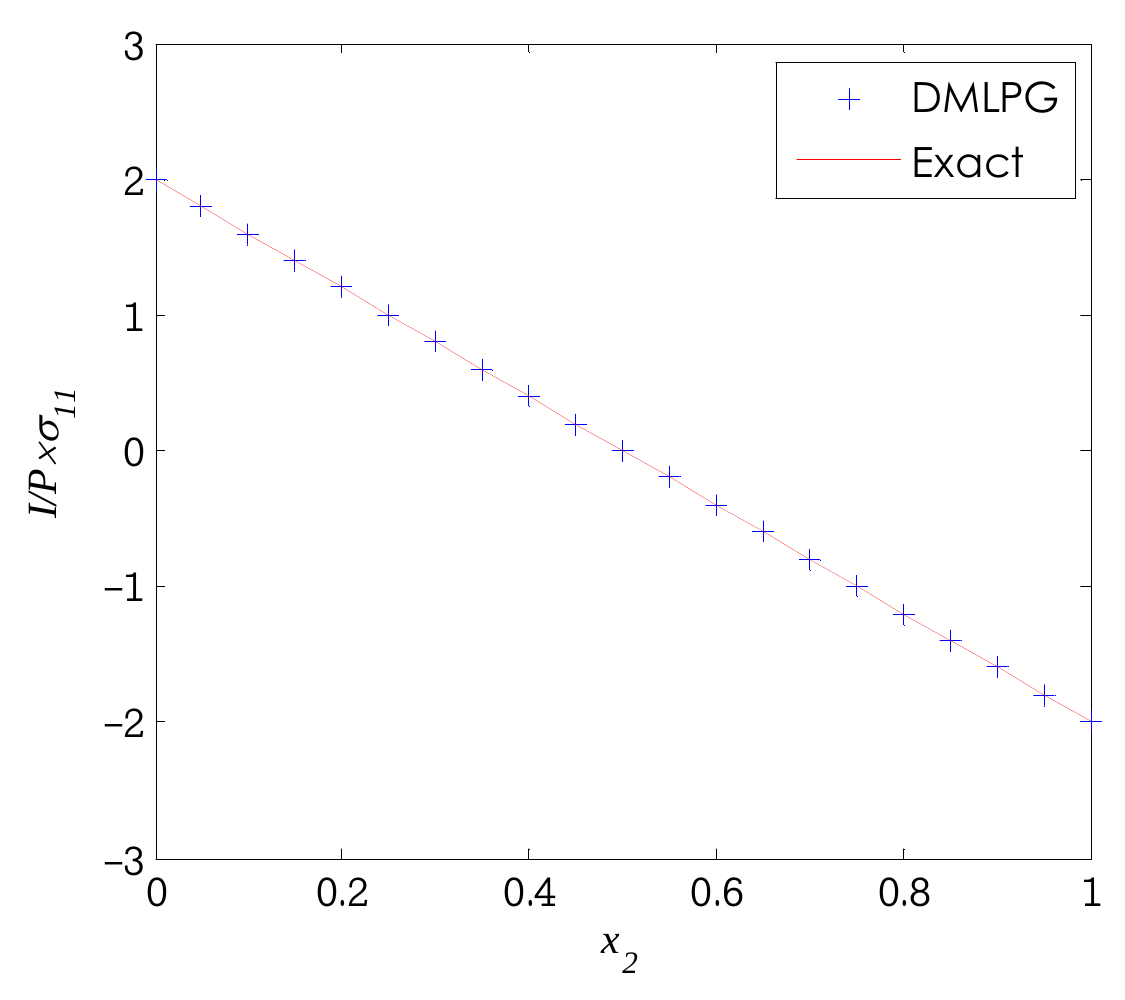}
\includegraphics[width=8cm]{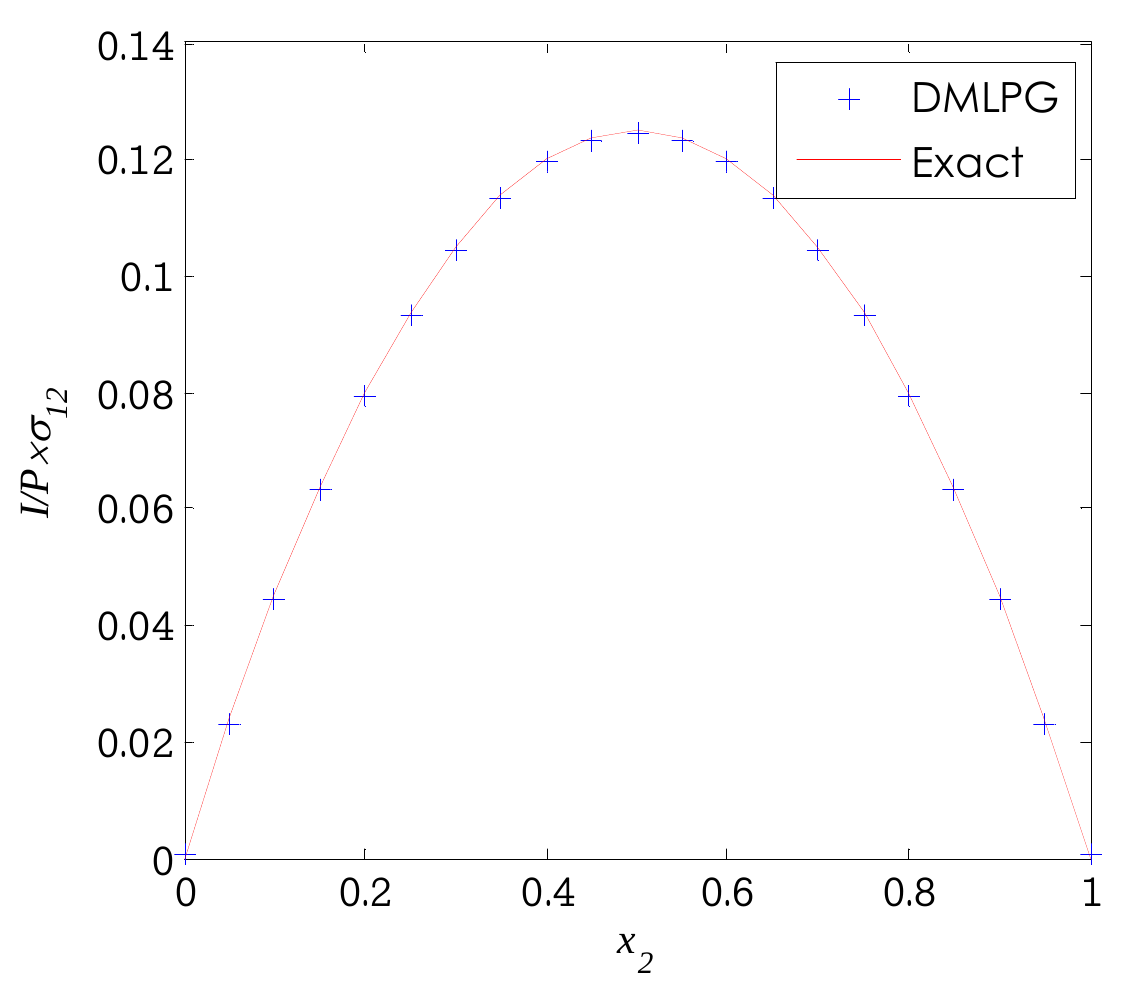}
\caption{\small{Numerical and exact normal and shear stresses at $x_1=4$ in cantilever beam}}\label{figure5}
\end{center}
\end{figure}

Finally, the DMLPG solutions of normal stress $\sigma_{11}$ and shear stress $\sigma_{12}$ at $\schi_1=L/2=4$ are plotted in Fig. \ref{figure5} and they are compared with the exact solutions.

\subsection{Infinite plate with circular hole}
\begin{figure}[hbt]
\begin{center}
\includegraphics[width=9cm]{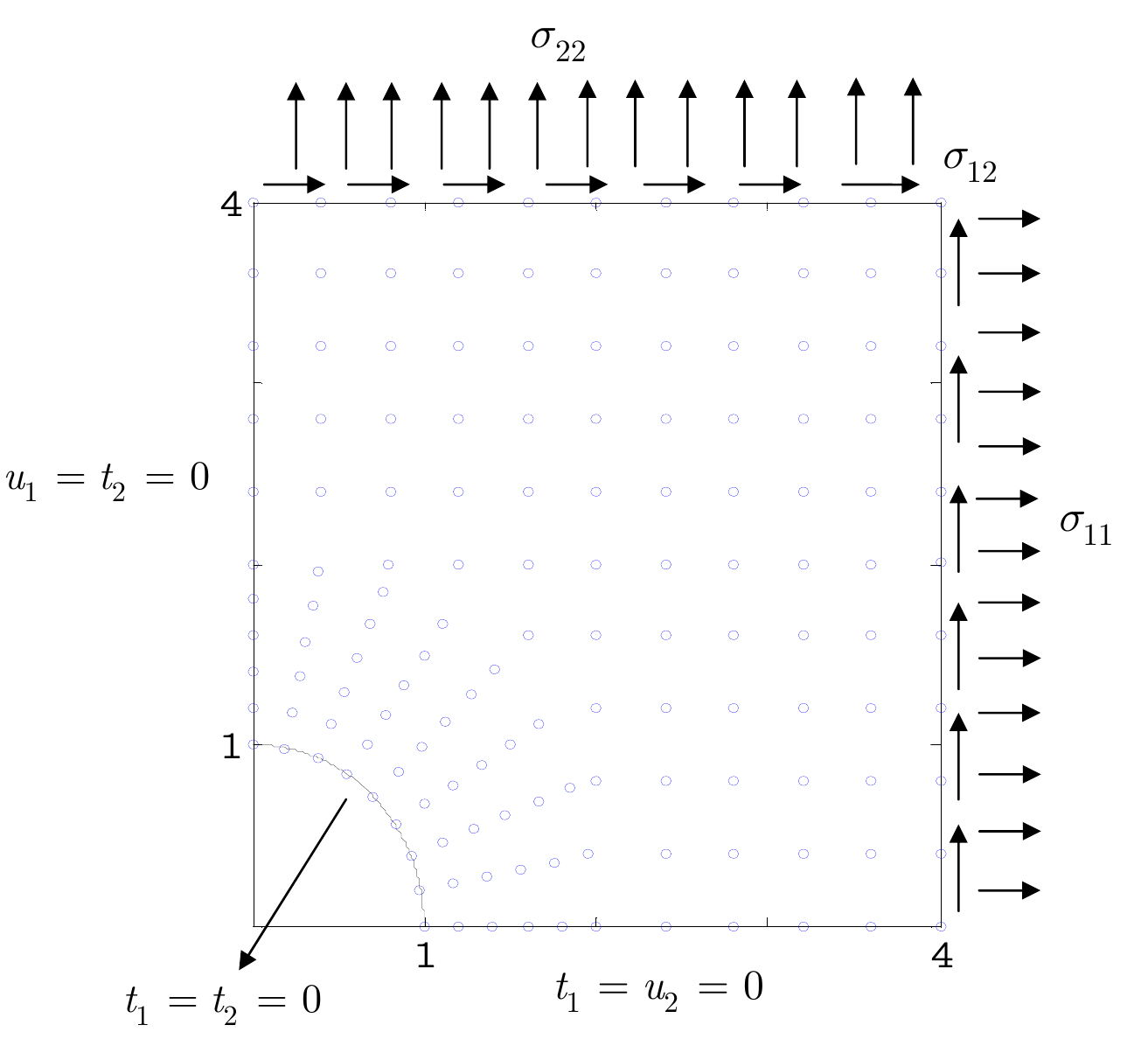}
\caption{\small{A quadrant of plate with circular hole, meshless points and boundary conditions}}\label{figure6}
\end{center}
\end{figure}

\begin{figure}[hbt]
\begin{center}
\includegraphics[width=8cm]{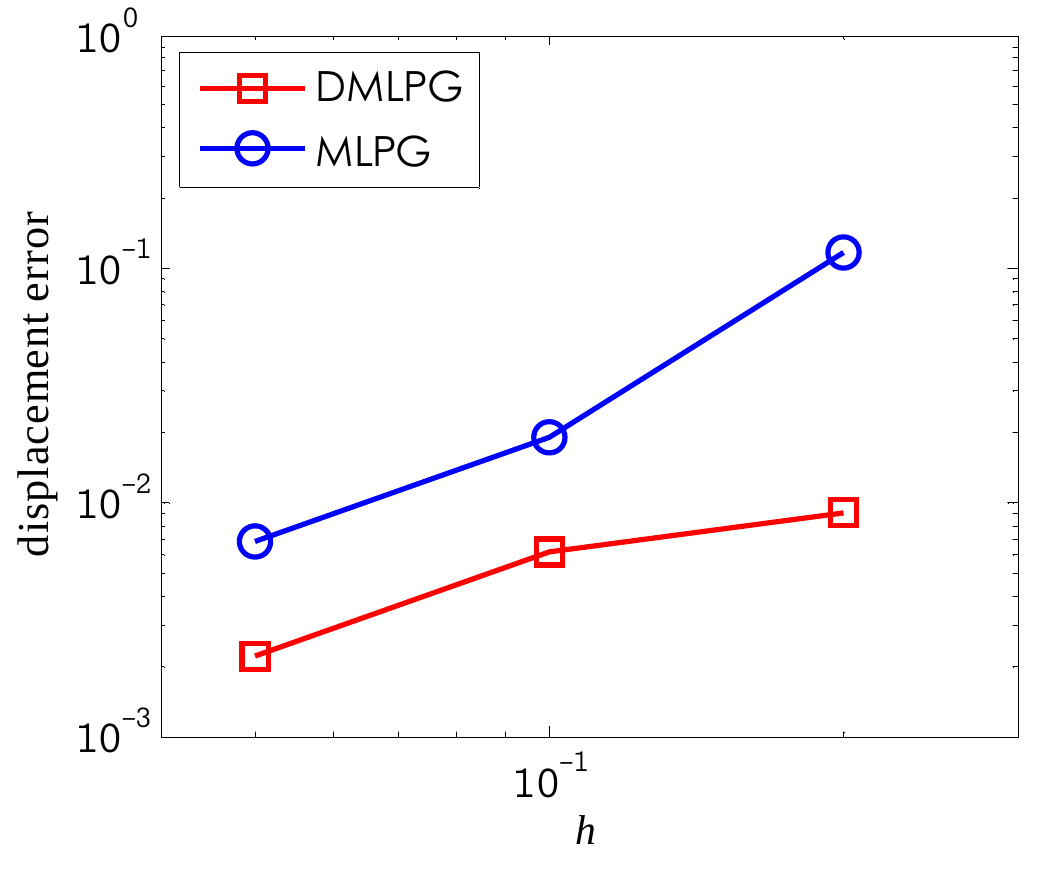}
\includegraphics[width=8cm]{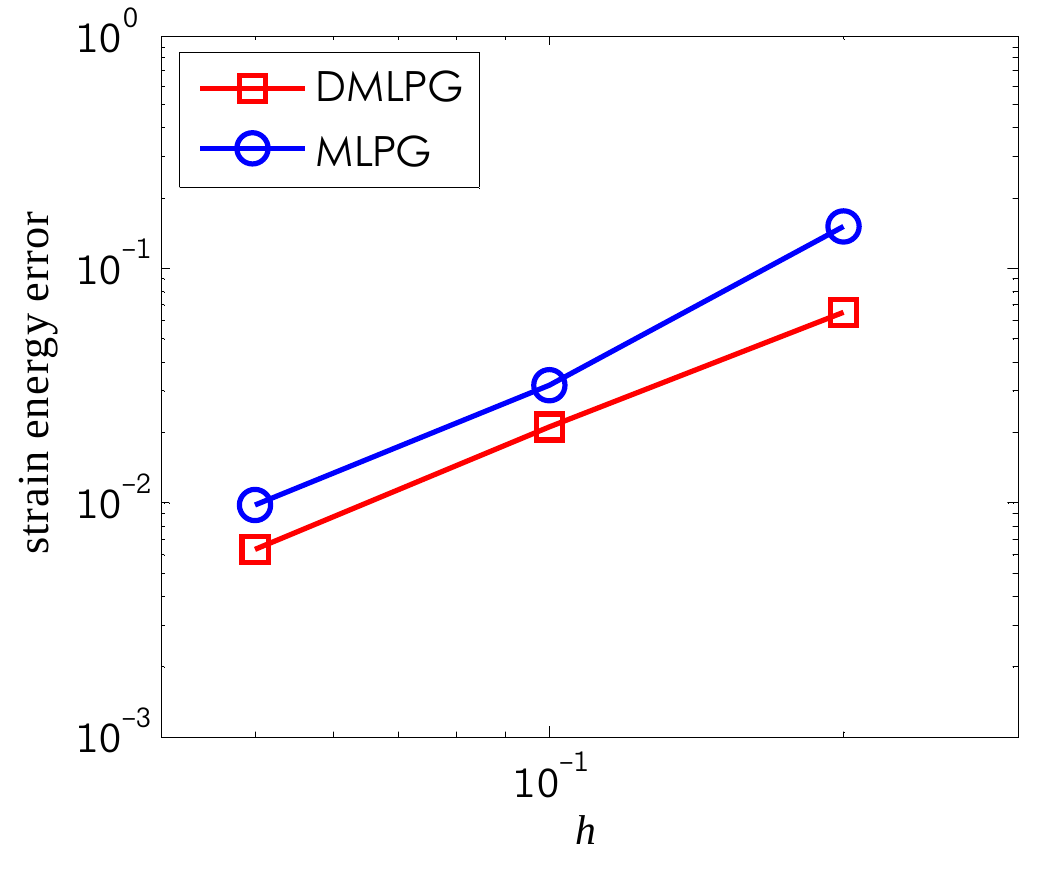}
\caption{\small{Relative displacement and strain errors for infinite plate with hole.}}\label{figure7}
\end{center}
\end{figure}

\begin{figure}[hbt]
\begin{center}
\includegraphics[width=8cm]{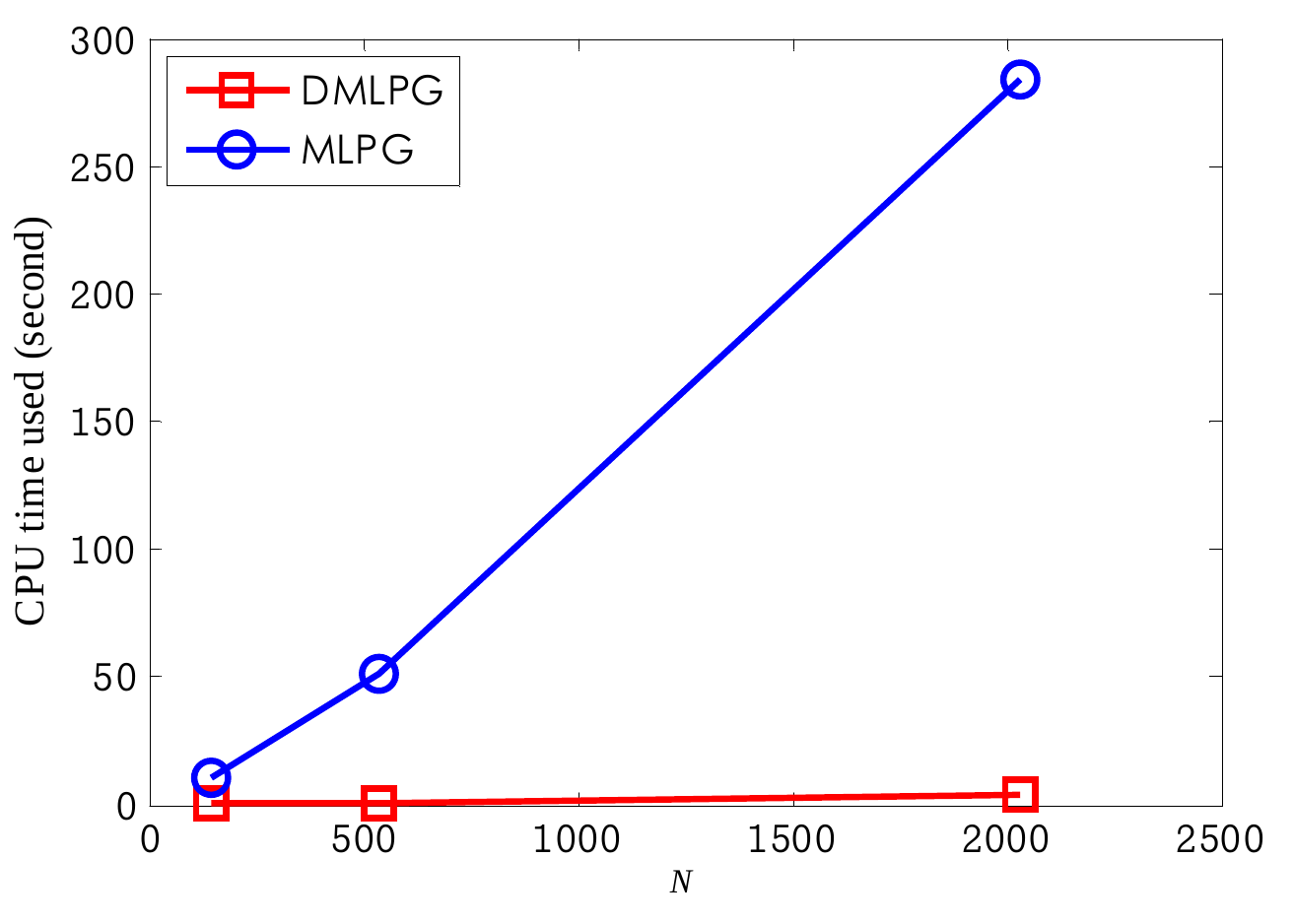}
\caption{\small{Computational costs for infinite plate with hole}}\label{figure8}
\end{center}
\end{figure}

\begin{figure}[hbt]
\begin{center}
\includegraphics[width=8cm]{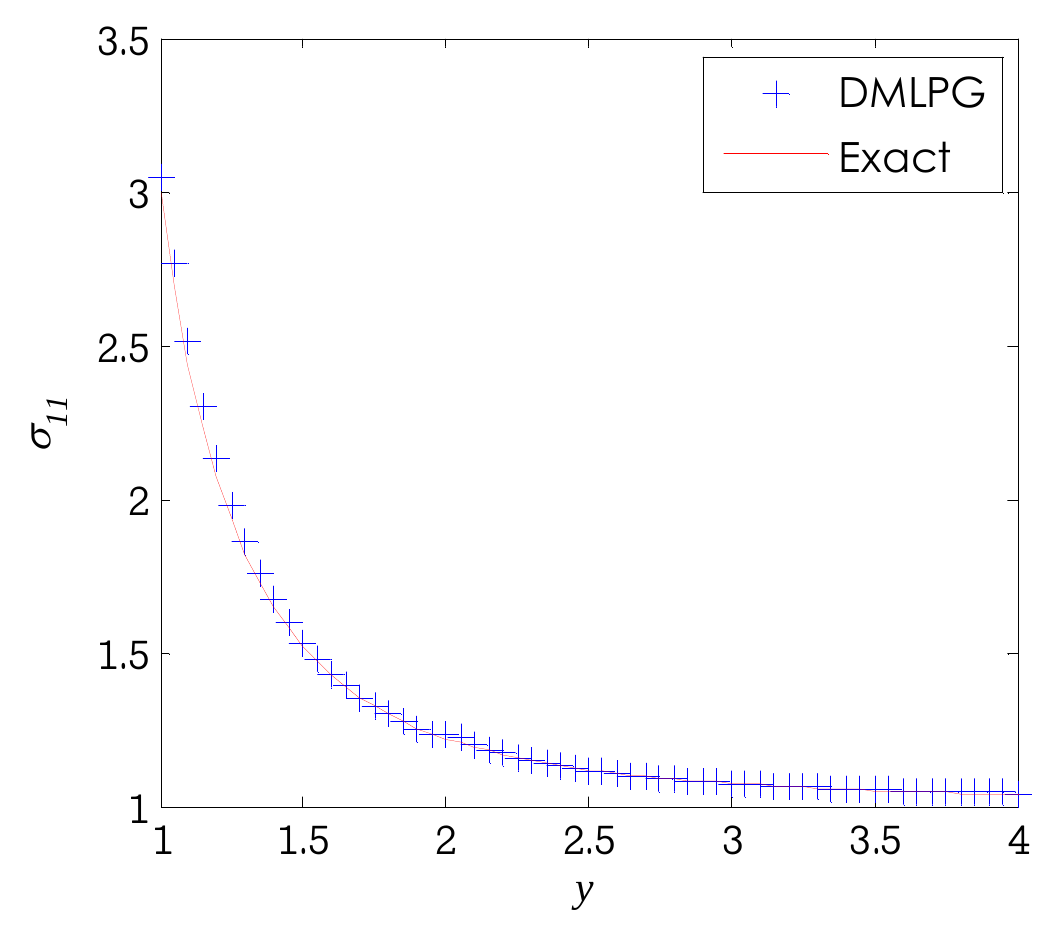}
\includegraphics[width=8cm]{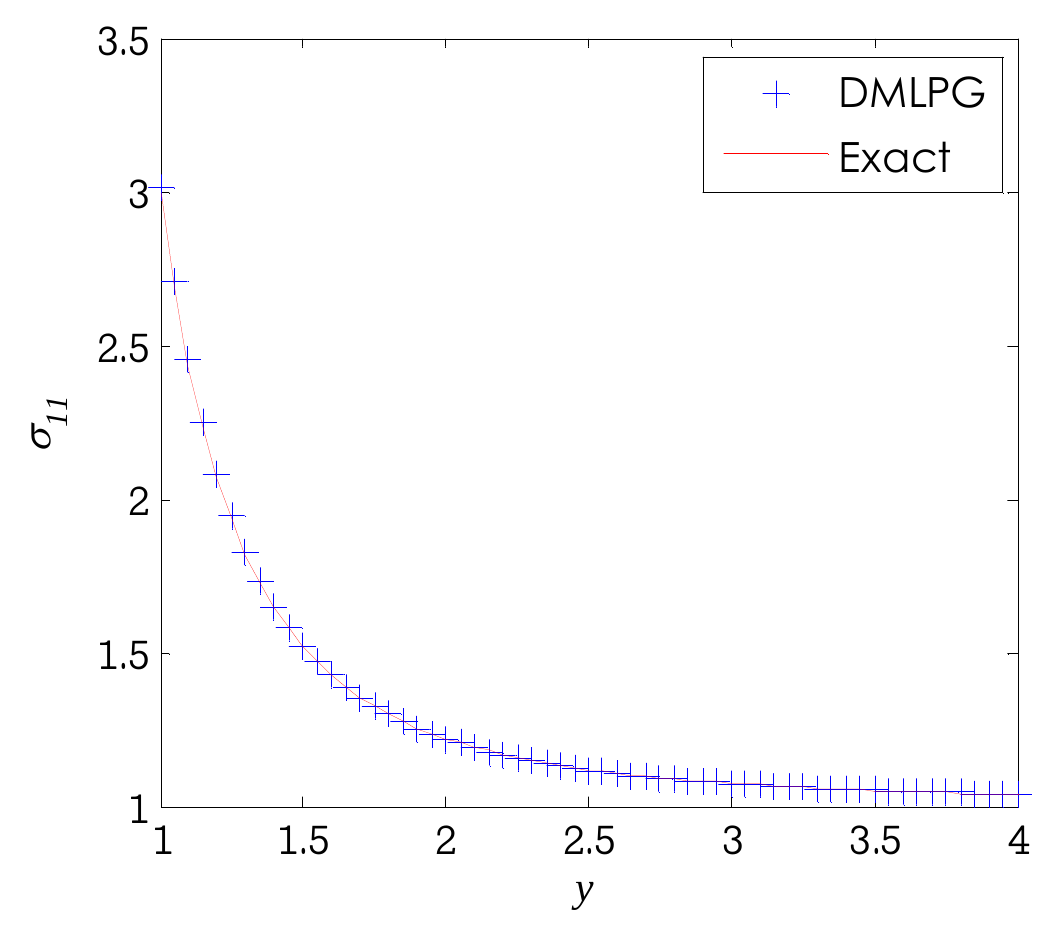}
\caption{\small{Numerical and exact normal stresses in plate, 535 nodes (up), 2034 nodes (down)}}\label{figure9}
\end{center}
\end{figure}

Consider an infinite plate with a central hole $\schi_1^2+\schi_2^2\leqslant a^2$ of
radius $a$, subjected to a unidirectional tensile load of $\sigma=1$ in the $\schi_1$-direction at infinity.
There is an analytical
solution for stress in the polar coordinate $(r,\theta)$
\begin{align*}
\sigma_{11}=& \sigma\left[ 1-\frac{a^2}{r^2}\left( \frac{3}{2}\cos 2\theta +\cos 4\theta\right)+\frac{3a^4}{2r^4}\cos 4\theta \right],\\
\sigma_{12}=& \sigma\left[ -\frac{a^2}{r^2}\left( \frac{1}{2}\sin 2\theta +\sin 4\theta\right)+\frac{3a^4}{2r^4}\sin 4\theta \right],\\
\sigma_{22}=& \sigma\left[ -\frac{a^2}{r^2}\left( \frac{1}{2}\cos 2\theta -\cos 4\theta\right)-\frac{3a^4}{2r^4}\cos 4\theta \right],
\end{align*}
with the corresponding displacements
\begin{align*}
u_{1}=& \frac{1+\overline \nu}{\overline E}\sigma\left[ \frac{1}{1+\overline \nu}r\cos\theta+\frac{2}{1+\overline \nu} \frac{a^2}{r}\cos \theta +
\frac{1}{2}\frac{a^2}{r}\cos 3\theta-\frac{1}{2}\frac{a^4}{r^3}\cos 3\theta \right],\\
u_{2}=& \frac{1+\overline \nu}{\overline E}\sigma\left[ \frac{-\nu}{1+\overline \nu}r\sin\theta-\frac{1-\nu}{1+\overline \nu} \frac{a^2}{r}\sin \theta +\frac{1}{2}\frac{a^2}{r}\sin 3\theta-\frac{1}{2}\frac{a^4}{r^3}\sin 3\theta \right].
\end{align*}
In computations, we consider a finite plate of length $b=4$ with a circular hole of radius $a=1$ (see Fig. \ref{figure6}), where the solution is very close to that of the
infinite plate \cite{roark-young:1975}.
Due to symmetry, only the upper right quadrant of the plate is
modelled.  The traction
boundary conditions given by the exact solution
are imposed on the right and top
edges (see Fig. \ref{figure6}).
Symmetry conditions are imposed on the left
and bottom edges, i.e., $u_1 = 0,\,t_2=0$ are prescribed on
the left edge and $u_2=0,\, t_1=0$ on the bottom edge, and the inner boundary at $a =1$ is traction free, i.e. $t_1=t_2=0$.
Numerical results are presented for a plane stress case with $E=1.0$ and $\nu=0.25$. The initial set point is depicted in Fig. \ref{figure6},
where we use more points near the hole. Thus the support size $\delta$ varies according to the density of neighboring points. Here $\delta=2mh$
and $\delta=2.5mh$ are used for points near the hole and points far away from the hole, respectively. Mesh-size $h$ is defined to be $\min\{h_r,h_\theta\}$ for the points close to the hole. In DMLPG, we use circular subdomains for points located on the arc boundary $r=a$,
and square subdomains for other points. Computations are repeated by halving $h_r$ and $h_\theta$, twice.
Results are presented in Figs. \ref{figure7} and \ref{figure8} which compare the displacement errors, the strain energy errors, and the CPU times used.
Moreover, the exact normal stress $\sigma_{11}$ at $\schi_1= 0$ is plotted in Fig. \ref{figure9} and it is compared with the DMLPG solution.

\subsection{3D Boussinesq problem}
The Boussinesq problem can be described as a concentrated load
acting on a semi-infinite elastic medium with no body
force. The exact displacement field within the semi-infinite
medium is given by Timoshenko and Goodier \cite{timoshenko-goodier:1970-1}
\begin{align*}
u_r &= \frac{(1+\nu)P}{2E \pi \rho}\left[\frac{zr}{\rho^2}-\frac{(1-2\nu)r}{\rho+z} \right],\\
w &= \frac{(1+\nu)P}{2E \pi \rho}\left[\frac{z^2}{\rho^2}+2(1-\nu) \right].
\end{align*}
where $u_r$ is the radial displacement, $w$ (or $u_3$) is the vertical displacement, $\rho=\sqrt{\schi_1^2+\schi_2^2+\schi_3^2}$ is the distance to the loading point and $r=\sqrt{\schi_1^2+\schi_2^2}$ is the
projection of $\rho$ on the loading surface. The exact stresses field is
\begin{align*}
\sigma_r &= \frac{P}{2 \pi \rho^2}\left[-\frac{3zr^2}{\rho^3}+\frac{(1-2\nu)\rho}{\rho+z} \right],\\
\sigma_\theta &=  \frac{(1-2\nu)P}{2 \pi \rho^2}\left[\frac{z}{\rho}-\frac{\rho}{\rho+z} \right],\\
\sigma_{zz}& =-\frac{3\pi z^3}{2\pi \rho^5} ,\\
\tau_{zr}&= \tau_{rz} = -\frac{3\pi rz^2}{2\pi \rho^5}.
\end{align*}
It is clear that the displacements and stresses are strongly
singular and they approach infinity; with the displacement
being $O(1/\rho)$
and the stresses being $O(1/\rho^2)$. MLPG has been applied to this problem in \cite{li-et-al:2003-1}.

\begin{figure}[hbt]
\begin{center}
\includegraphics[width=8cm]{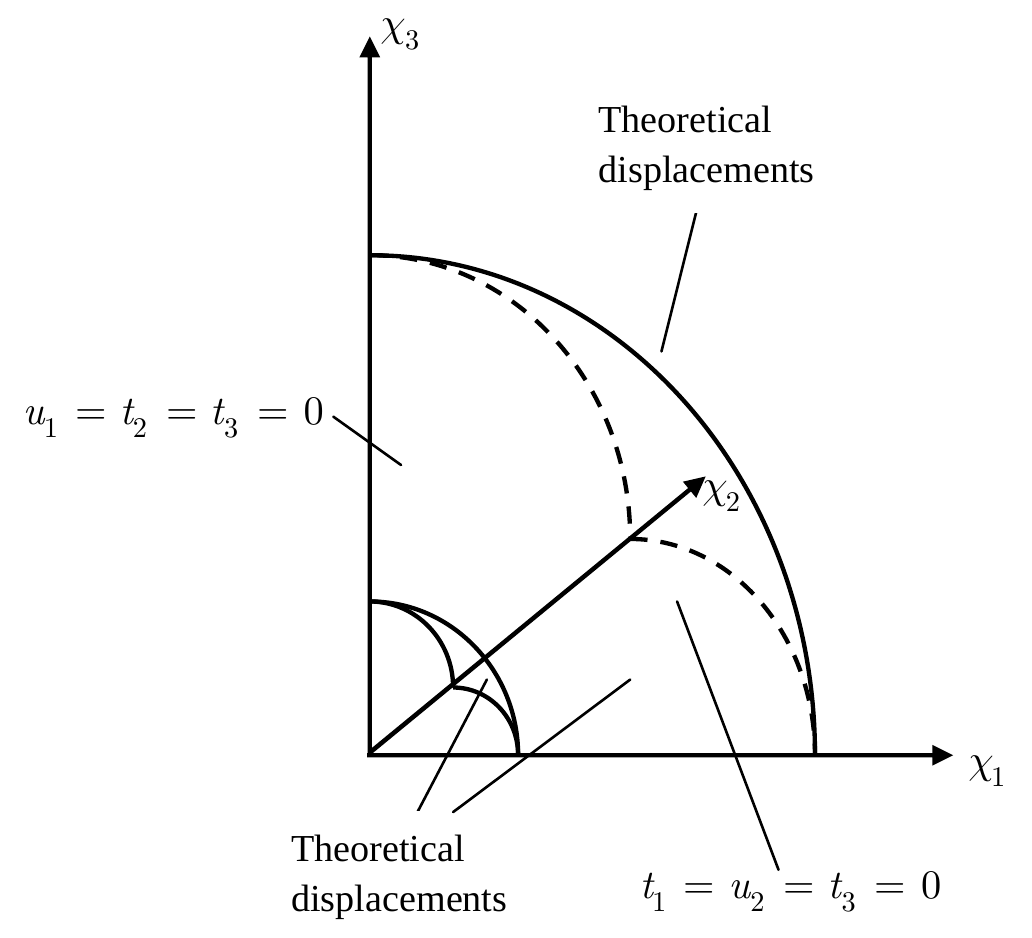}
\includegraphics[width=8cm]{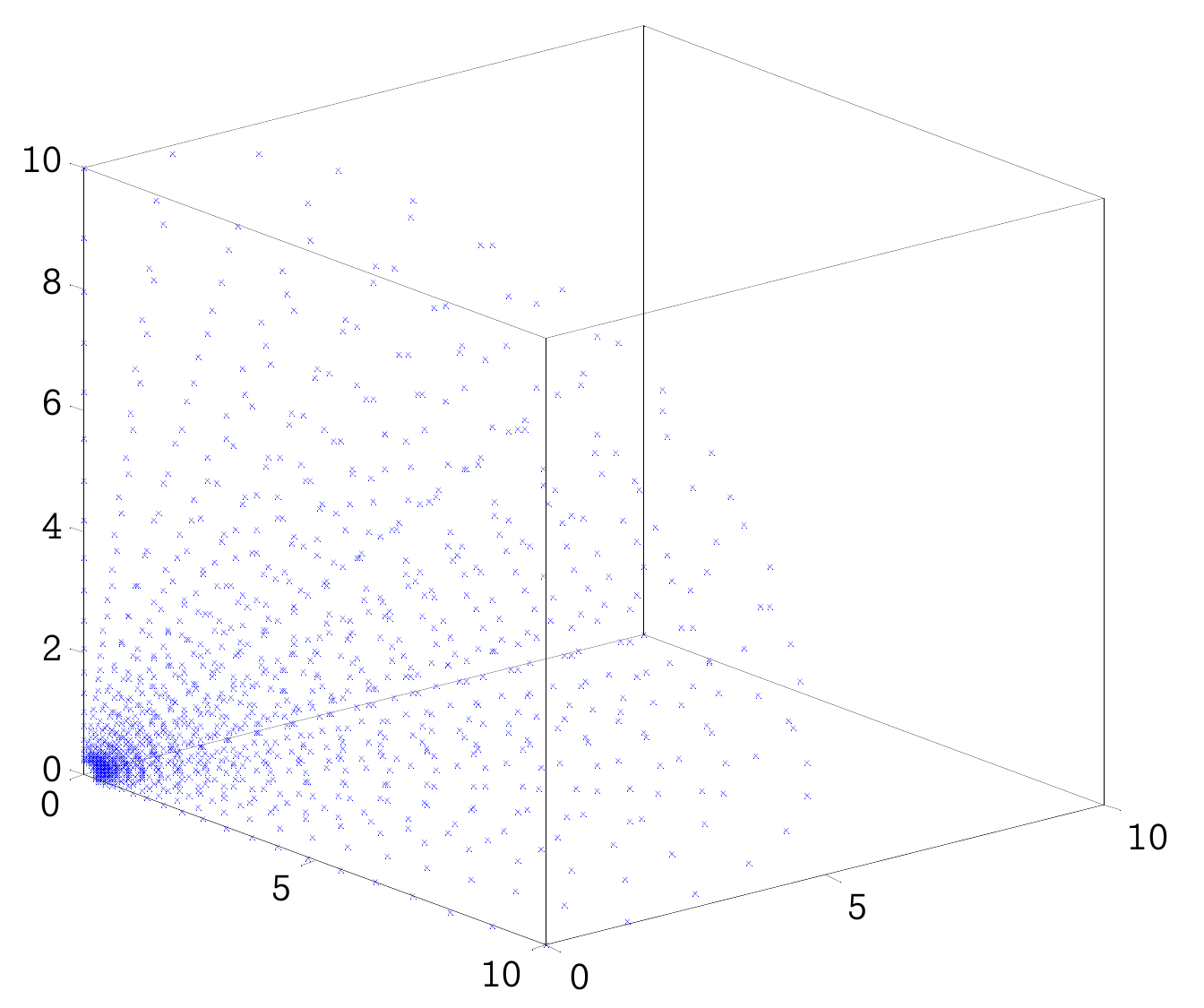}
\caption{\small{The consideration domain and meshless points (1386 points) in Boussinesq problem}}\label{figure10}
\end{center}
\end{figure}

\begin{figure}[hbt]
\begin{center}
\includegraphics[width=8cm]{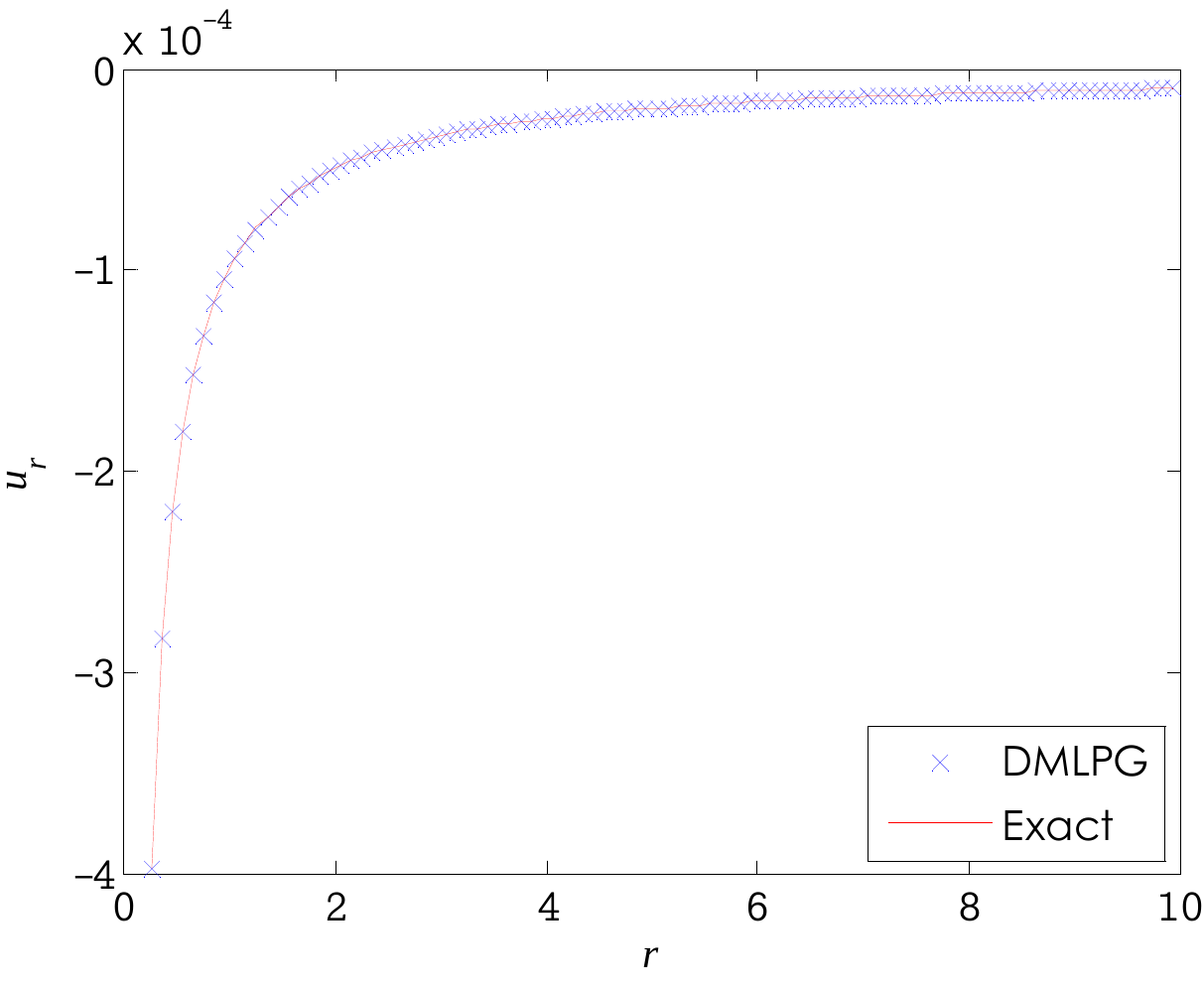}
\includegraphics[width=8cm]{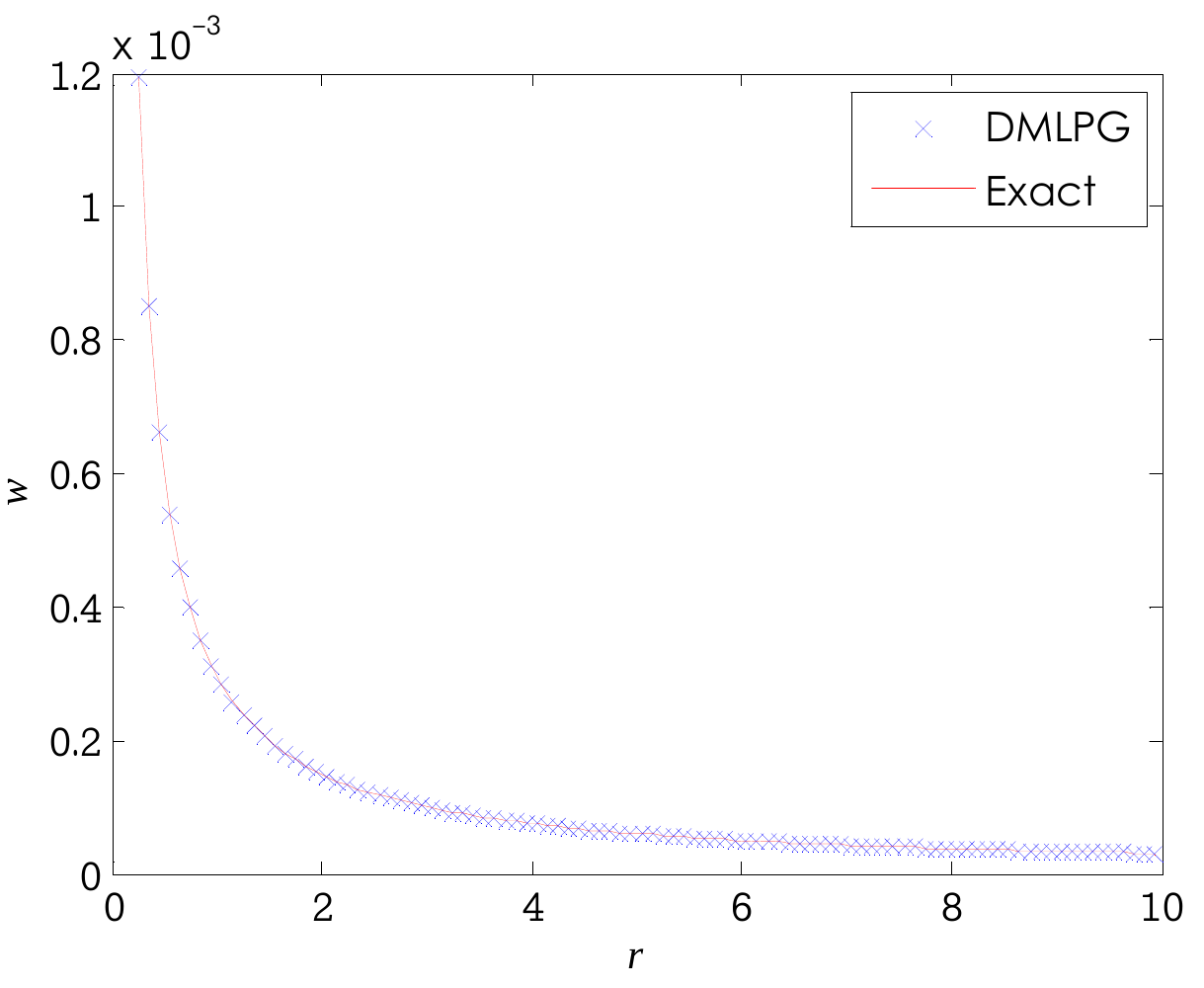}
\caption{\small{Radial Displacement $u_r$ and vertical displacement $w$
in loading surface in Boussinesq problem}}\label{figure11}
\end{center}
\end{figure}

\begin{figure}[hbt]
\begin{center}
\includegraphics[width=8cm]{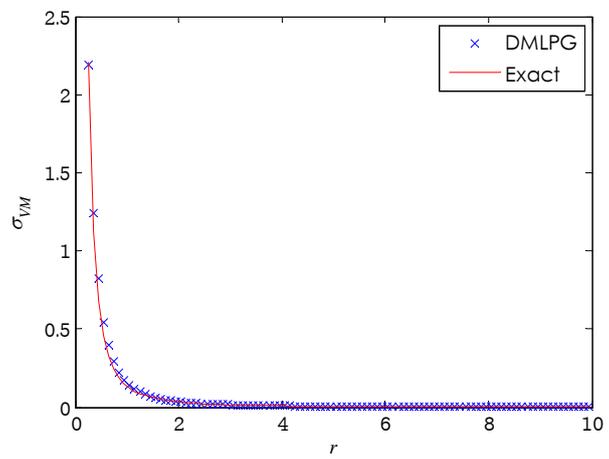}
\caption{\small{Von Mises Stress in loading surface in Boussinesq problem}}\label{figure12}
\end{center}
\end{figure}

In numerical simulation, a finite sphere with large radius $b=10$ is used. Due to the symmetry, a first one-eighth of the sphere is considered
and symmetry boundary conditions are applied on planes $xz$ and $yz$ (see Fig. \ref{figure10}). In fact we impose $t_1=u_2=t_3=0$ on plane $xz$, and
$u_1=t_2=t_3=0$ on plane $yz$. In order to avoid direct encounter with
the singular loading point, the theoretical displacement
is applied on a small spherical surface with radius $b/40=0.25$. An isotropic material of $E = 1000$, $\nu = 0.25$ and $P=1$ is used.
The number of meshless points is $1386$, which are scattered inside the domain and on the boundary. The density of nodes depends on the distance from the loading points, where we have many points near the small sphere and few points far from it (see Fig. \ref{figure10}).
Thus the support size $\delta$ varies and depends on $\rho$, correspondingly.
Analytical and DMLPG solutions of the radial displacement $u_r$ and vertical displacement $w$
on the surface $xy$ are plotted
in Fig. \ref{figure11}. The Von Mises stress on the surface $xy$ is also shown in
Fig. \ref{figure12}. These are the results of DMLPG1 with cubes as sub-domains where the CPU time used is around 3 seconds.
Again we note that a 2-point Gaussian quadrature in each axis gives the exact numerical integration.
The same results will be obtained by DMLPG5.

Finally for comparison we apply both MLPG1 and MLPG5 to this problem with the same meshless points and MLS parameters.
The accuracy of results are far less than DMLPG solutions and the CPU run times are about 7400 sec. for MLPG1 and 450 sec. for MLPG5.
In computations, a 10-point Gauss formula is employed in each axis. In fact, for MLPG1, the MLS shape function subroutines should be called 1000 times
to integrate over a sub-domain $\Omega_k$. In MLPG5 this number reduces to 100, because the integrals are all boundary integrals in this example.
Compare with DMLPG where the MLS subroutines are not called for integrations at all, leading to 3 sec. running time in this example.


\section{Conclusion}
In this paper we developed a new meshfree method for elasticity problems, which is
a weak form method in the cost-level of collocation (integration-free) methods. Integrations have been shifted into the MLS itself, rather than into an outside loop over
calls to MLS routines. In fact, we need to integrate against low-degree polynomials basis functions instead of complicated MLS shape functions. Besides, in some situations we can perform {\em exact} numerical integrations.
We applied DMLPG1 and 5 for
problems in two and three dimensional elasticity in this paper. The new methods can be easily applied to other problems in solid engineering.
On a downside, DMLPG1 and 5 do not work for linear basis functions ($m=1$). In addition, because of symmetry properties of polynomials in local subdomains, \cite{mirzaei-schaback:2013-1} shows that the convergence rates do not increase when going from $m=2k$ to $m=2k+1$. But the results show that this observation affects MLPG and DMLPG in the same way.
DMLPG4 can be formulated using the strategy presented in \cite{atluri-et-al:2000-1} to make the second unsymmetric local weak forms and applying the GMLS approximation of this paper.
Finally, we believe that DMLPG methods have great potential to replace the original
MLPG methods in many situations, specially for three dimensional problems.

\section*{Acknowledgment}
Special thanks go to Prof. R. Schaback, Universit\"at G\"ottingen, Dr. K. Hasanpour, Department of Mechanical Engineering , University of Isfahan, and Dr. K. Mohajer for their useful helps and comments.


\begin{thebibliography}{10}
\expandafter\ifx\csname url\endcsname\relax
  \def\url#1{\texttt{#1}}\fi
\expandafter\ifx\csname urlprefix\endcsname\relax\def\urlprefix{URL }\fi
\expandafter\ifx\csname href\endcsname\relax
  \def\href#1#2{#2} \def\path#1{#1}\fi

\bibitem{atluri-zhu:1998-1}
S.~Atluri, T.-L. Zhu, A new meshless local {P}etrov-{G}alerkin ({MLPG})
  approach in computational mechanics, Computational Mechanics 22 (1998)
  117--127.

\bibitem{atluri-zhu:2000-1}
S.~N. Atluri, T.~L. Zhu, The meshless local {P}etrov-{G}alerkin ({MLPG})
  approach for solving problems in elasto-statics, Computational Mechanics 25
  (2000) 169--179.

\bibitem{li-et-al:2003-1}
Q.~Li, S.~Shen, Z.~D. Han, S.~N. Atluri, Application of meshless local
  {P}etrov-{G}alerkin ({MLPG}) to problems with singularities, and material
  discontinuities, in 3-{D} elasticity, CMES: Computer Modeling in Engineering
  \& Sciences 4 (2003) 571--585.

\bibitem{sladek-et-al:2005-1}
J.~Sladek, V.~Sladek, C.~Zhang, An advanced numerical method for computing
  elastodynamic fracture parameters in functionally graded materials,
  Computational Materials Science 32 (2005) 532--543.

\bibitem{sladek-et-al:2013-1}
J.~Sladek, P.~Stanak, Z.~D. Han, V.~Sladek, S.~N. Atluri, {A}pplications of the
  {MLPG} method in engineering \& sciences: {A} review, {CMES}--Computer
  Modeling in Engineering \& Sciences 92 (2013) 423--475.

\bibitem{belytschko-et-al:1994-1}
T.~Belytschko, Y.~Lu, L.~Gu, Element-{F}ree {G}alerkin methods, International
  Journal for Numerical Methods in Engineering 37 (1994) 229--256.

\bibitem{beissel-belytsschko:1996-1}
S.~Beissel, T.~Belytschko, Nodal integration of the element-free {G}alerkin
  method, Computer Methods in Applied Mechanics and Engineering 139 (1996)
  49--74.

\bibitem{dolbow-belytsschko:1999-1}
J.~Dolbow, T.~Belytschko, Numerical integration of the {G}alerkin weak form in
  meshfree methods, Computational Mechanics 23 (1999) 219--230.

\bibitem{atluri-et-al:1999-2}
S.~N. Atluri, H.~G. Kim, J.~Y. Cho, A critical assessment of the truly
  {M}eshless {L}ocal {P}etrov-{G}alerkin ({MLPG}), and {L}ocal {B}oundary
  {I}ntegral {E}quation ({LBIE}) methods, Computational Mechanics 24 (1999)
  348--372.

\bibitem{carpinteri-et-al:2002-1}
A.~Carpinteri, G.~Ferro, G.~Ventura, The partition of unity quadrature in
  meshless methods, International Journal for Numerical Methods in Engineering
  54 (2002) 987--1006.

\bibitem{pecher:2006-1}
R.~Pecher, Efficient cubature formulae for {MLPG} and related methods,
  International Journal for Numerical Methods in Engineering 65 (2006)
  566--593.

\bibitem{mazzia-pini:2010-1}
A.~Mazzia, G.~Pini, Product {G}auss quadrature rules vs. cubature rules in the
  meshless local {P}etrov-{G}alerkin method, Journal of Complexity 26 (2010)
  82--101.

\bibitem{babuska-et-al:2009-1}
I.~Babuska, U.~Banerjee, J.~Osborn, Q.~Zhang, Effect of numerical integration
  on meshless methods, Comput. Methods Appl. Mech. Engrg. 198 (2009) 27--40.

\bibitem{mirzaei-schaback:2013-1}
D.~Mirzaei, R.~Schaback, Direct {M}eshless {L}ocal {P}etrov-{G}alerkin
  ({DMLPG}) method: a generalized {MLS} approximation, Applied Numerical
  Mathematics 33 (2013) 73--82.

\bibitem{mirzaei-et-al:2012-1}
D.~Mirzaei, R.~Schaback, M.~Dehghan, On generalized moving least squares and
  diffuse derivatives, IMA Journal of Numerical Analysis 32 (2012) 983--1000.

\bibitem{mirzaei:2013-1}
D.~Mirzaei, Error boounds for {GMLS} derivatives approximations of {S}obolev
  functions, preprint, University of Isfahan, Available at
  \verb+http://sci.ui.ac.ir/+$\sim$\verb+d.mirzaei+ (2014).

\bibitem{mirzaei-schaback:2013-2}
D.~Mirzaei, R.~Schaback, Solving heat conduction problem by the {D}irect
  {M}eshless {L}ocal {P}etrov-{G}alerkin ({DMLPG}) method, Numerical Algorithms
  65 (2014) 275--291.

\bibitem{mazzia-et-al:2012-1}
A.~Mazzia, G.~Pini, F.~Sartoretto, Numerical investigation on direct {MLPG} for
  {2D} and {3D} potential problems, CMES: Computer Modeling in Engineering \&
  Sciences 88 (2012) 183--209.

\bibitem{wendland:2005-1}
H.~Wendland, Scattered Data Approximation, Cambridge University Press, 2005.

\bibitem{lancaster-salkauskas:1981-1}
P.~Lancaster, K.~Salkauskas, Surfaces generated by moving least squares
  methods, Mathematics of Computation 37 (1981) 141--158.

\bibitem{belytschko-et-al:1996-1}
T.~Belytschko, Y.~Krongauz, D.~Organ, M.~Fleming, P.~Krysl, Meshless methods:
  an overview and recent developments, Computer Methods in Applied Mechanics
  and Engineering, special issue 139 (1996) 3--47.

\bibitem{mirzaei:2015-1}
D.~Mirzaei, Analysis of moving least squares approximation revisited, Journal
  of Computational and Applied Mathematics (2015) In press.

\bibitem{atluri-et-al:1999-1}
S.~N. Atluri, J.~Y. Cho, H.~G. Kim, Analysis of thin beams, using the meshless
  local {P}etrov-{G}alerkin method, with generalized moving least squares
  interpolations, Computational Mechanics 24 (1999) 334--347.

\bibitem{timoshenko-goodier:1970-1}
S.~P. Timoshenko, J.~N. Goodier, Theory of Elasticity, 3rd edition,
  McGraw-Hill, New York, 1970.

\bibitem{roark-young:1975}
R.~J. Roark, W.~C. Young, Formulas for Stress and Strain, McGraw-Hill, 1975.

\bibitem{atluri-et-al:2000-1}
S.~N. Atluri, J.~Sladek, V.~Sladek, T.-L. Zhu, The local boundary integral
  equation ({LBIE}) and it's meshless implementation for linear elasticity,
  Computational Mechanics 25 (2000) 180--198.

\end{thebibliography}

\newpage

\end{document}